# REPRESENTATION OF SOLUTIONS TO BSDES ASSOCIATED WITH A DEGENERATE FSDE

By Jianfeng Zhang

*University of Southern California*

In this paper we investigate a class of decoupled forward–backward SDEs, where the volatility of the FSDE is *degenerate* and the terminal value of the BSDE is a *discontinuous* function of the FSDE. Such an FBSDE is associated with a degenerate parabolic PDE with discontinuous terminal condition. We first establish a Feynman–Kac type representation formula for the spatial derivative of the solution to the PDE. As a consequence, we show that there exists a stopping time $\tau$ such that the martingale integrand of the BSDE is continuous before $\tau$ and vanishes after $\tau$. However, it may blow up at $\tau$, as illustrated by an example. Moreover, some estimates for the martingale integrand before $\tau$ are obtained. These results are potentially useful for pricing and hedging discontinuous exotic options (e.g., digital options) when the underlying asset's volatility is small, and they are also useful for studying the rate of convergence of finite-difference approximations for degenerate parabolic PDEs.

**1. Introduction.** In this paper we investigate the following decoupled forward–backward SDE:

$$(1.1) \quad \begin{aligned} X_t &= x + \int_0^t b(r, X_r)\, dr + \int_0^t \sigma(r, X_r)\, dW_r; \\ Y_t &= g(X_T) + \int_t^T f(r, X_r, Y_r, Z_r)\, dr - \int_t^T Z_r\, dW_r, \end{aligned}$$

where $W$ is a standard Brownian motion, $\sigma, b, f$ and $g$ are deterministic functions. It is well known that, in mathematical finance theory, the solution triple $(X, Y, Z)$ can be interpreted as underlying asset price, option price and hedging strategy, respectively (see, e.g., [2]). The equations of type (1.1)

Received July 2004; revised November 2004.
Supported in part by NSF Grant DMS-04-03575.

were first studied by Pardoux and Peng [9]. We refer the readers to the book of Ma and Yong [5] for more details on the subject. Among other things, Pardoux and Peng [9] showed that (1.1) was related to the following quasilinear parabolic PDE:

$$u_t + \tfrac{1}{2}\sigma^2 u_{xx} + b u_x + f(t,x,u,u_x\sigma) = 0;$$
(1.2)
$$u(T,x) = g(x),$$

in the sense that

(1.3) $$Y_t = u(t, X_t), \qquad Z_t = (u_x\sigma)(t, X_t).$$

For the purpose of applications, we are particularly interested in pathwise properties of the process $Z$. In the literature there are typically two types of conditions to ensure the regularity of $Z$. One is to assume that the coefficients $b, \sigma, f$ and $g$ are sufficiently smooth (e.g., [9] and [4]) so that (1.2) has a classical solution $u$ and, thus, $Z$ is continuous. The other is to assume that $\sigma$ is uniformly nondegenerate (e.g., [3], [6] and [7]) so that $X_T$ has a density (see, e.g., [8]) and, thus, $u(t,x)$ is smooth in $x$ for $t < T$, thanks to the nonlinear Feynman–Kac formula.

It is our goal of this paper to remove both conditions above. We will allow $\sigma$ to be *degenerate* and $g$ to be *discontinuous*. We note that in this case (1.2) is a degenerate PDE which, in general, has no smooth solution. A trivial counterexample is that $\sigma = b = f = 0$ and $g = \mathbb{1}_{\{x>0\}}$, then $u(t,x) = \mathbb{1}_{\{x>0\}}$ for $\forall t \in [0,T]$, and, hence, $u$ is discontinuous in $x$. However, by (1.3) and noting that $\sigma = 0$, one may still view $Z_t = 0$ in this example. In fact, $Z_t = 0$ is indeed the solution to (1.1).

A less obvious example is Example 1 in Section 4.1, in which the process $Z$ blows up at and only at some time $t < T$. It turns out that this is already the worst case one might encounter. We will show that, under certain conditions, there exists a stopping time $\tau$ such that $Z_t$ is continuous for $t < \tau$ and $Z_t = 0$ for $t \geq \tau$. So along each path, the only possible discontinuous point of $Z_t$ is $\tau$. Moreover, we have an explicit rule to locate $\tau$ and we have an estimate for $Z_t$ when $t < \tau$.

The main tool of our approach is a new Feynman–Kac type representation formula for $u_x$ (and $Z$) by using Malliavin calculus. As in [3] and [6], this formula does not involve the derivatives of $f$ or $g$ (thus, $g$ can be discontinuous!). But unlike those two works which require $\sigma$ to be uniformly nondegenerate, our new formula allows $\sigma$ to be degenerate. As a payoff, due to this degeneracy, our estimates for $u_x$ are technically much more involved than those in [6].

At this point we would like to mention that the discontinuity of $g$ is mainly motivated by digital options for which $g(x) = \mathbb{1}_{\{x>K\}}$. Degenerate diffusion also appears quite often in applications (noting that even in the standard

Black–Scholes model, the stock price equation is degenerate!). For example, in option pricing theory one may face a situation where the underlying asset market is quite stable during a random time period (so $\sigma$ is small) or there is a risk that the underlying corporation may go bankrupt at some random time (so $\sigma = 0$ afterward). Our results are potentially useful for pricing and hedging options in these markets. Also, the regularity of $u$ plays a very important role for studying the rate of convergence of finite-difference approximations for degenerate PDE (1.2) (see, e.g., [11]).

For technical reasons, in this paper we assume that all processes are one-dimensional and that $f$ is linear on $Z$. More general cases are left for future research. We note that in a recent paper Bally [1] studied the density of a degenerate multidimensional diffusion. We hope that his work may bring us some insights into our problem.

The rest of the paper is organized as follows. In Section 2 we give all the necessary preparations. In Section 3 we study two "good" cases which extend some results of [9] and [6], respectively. In Section 4 we study the case that $f = 0$ and derive a new representation formula for $Z$. Finally, in Section 5 we study the general case.

**2. Preliminaries.** Let $(\Omega, \mathcal{F}, P)$ be a complete probability space on which is defined a one-dimensional Brownian motion $W = (W_t)_{t \geq 0}$, and $\mathbf{F} \stackrel{\triangle}{=} \{\mathcal{F}_t\}_{t \geq 0}$ be the natural filtration generated by $W$, augmented by the $P$-null sets of $\mathcal{F}$.

The following spaces will be frequently used in the sequel: Let $O$ be an open subset of $[0, T] \times \mathbb{R}^k$ for some integer $k$,

- $C(O)$ is the space of all Lebesgue measurable functions on $[0, T] \times \mathbb{R}^k$ such that they are continuous in $O$;
- $C^{0,1}(O)$ is the space of those $\varphi \in C(O)$ such that they are continuously differentiable on the spatial variable(s) in $O$;
- $C_b^{0,1}(O)$ is the space of those $\varphi \in C^{0,1}(O)$ such that all the partial derivatives in $O$ are uniformly bounded (but $\varphi$ itself can be unbounded).

When $O = [0, T] \times \mathbb{R}^k$, we omit it. For example, $C^{0,1} = C^{0,1}([0, T] \times \mathbb{R}^k)$.

In this paper we assume all the processes involved are one-dimensional; and we shall use the following *Standing Assumptions*:

(A1) $b, \sigma \in C_b^{0,1}$;
(A2) $\sigma$ is uniformly continuous in $t$;
(A3) $f \in C([0, T] \times \mathbb{R}^3)$, and $f$ is uniformly Lipschitz continuous in $x, y, z$;
(A4) $g$ is Lebesgue measurable and $|g(x)| \leq \psi(x) \stackrel{\triangle}{=} K(1 + |x|^{p_0})$ for some constant $K$ and some $p_0 \geq 1$.

We note that, by assuming (A1),

(2.1) $$E\left\{\sup_{0 \leq t \leq T} \psi(X_t)^p\right\} \leq C_p \psi(x)^p \qquad \forall p \geq 1.$$

In fact, this is the only property of $\psi$ we will utilize in the rest of the paper. We also note that (A2) is equivalent to $\lim_{\varepsilon \to 0} \delta(\varepsilon) = 0$, where

$$(2.2) \quad \delta(\varepsilon) \triangleq \inf\left\{|t_1 - t_2| : 0 \leq t_1, t_2 \leq T, \sup_{x \in \mathbb{R}} |\sigma(t_1, x) - \sigma(t_2, x)| > \varepsilon\right\}.$$

Obviously, for any $x \in \mathbb{R}$ and $|t_1 - t_2| \leq \delta(\varepsilon)$, we have $|\sigma(t_1, x) - \sigma(t_2, x)| \leq \varepsilon$.

In order to simplify the presentation, we will also adopt the following assumption:

(A5) $\sigma, b$ are bounded.

However, without assuming it, all the results in the paper still hold true after some slight modification (see Remark 5.4).

Throughout the paper, we use a generic constant $K$ to denote all the Lipschitz constants involved. We also assume that $|\sigma(t,x)| + |b(t,x)| + |f(t,0,0,0)| \leq K$. Moreover, we use positive constants $C$ and $c$, which may vary from line to line but depend only on $K, T$ and the function $\psi$ in (A4), to denote upper bounds and lower bounds of estimates, respectively. Furthermore, if the bounds depend on some $p$ as well, we denote them by $C_p$ and $c_p$, respectively.

We now review some basic results, especially those concerning $Z$, in the literature. First, for any $(t,x) \in [0,T] \times \mathbb{R}$, let $(X_s^{t,x}, Y_s^{t,x}, Z_s^{t,x})_{t \leq s \leq T}$ denote the solution to the following FBSDE:

$$(2.3) \quad \begin{aligned} X_s^{t,x} &= x + \int_t^s b(r, X_r^{t,x})\,dr + \int_t^s \sigma(r, X_r^{t,x})\,dW_r; \\ Y_s^{t,x} &= g(X_T^{t,x}) + \int_s^T f(r, X_r^{t,x}, Y_r^{t,x}, Z_r^{t,x})\,dr - \int_s^T Z_r^{t,x}\,dW_r. \end{aligned}$$

When $t = 0$, (2.3) is the same as (1.1), and we still use $(X, Y, Z)$ to denote its solution. Next, we define $u(t,x) \triangleq Y_t^{t,x}$. It is well known that, under certain conditions, $u$ is the unique viscosity solution to (1.2) and $Y_t = u(t, X_t)$. Moreover, if $u \in C^{0,1}$, then (1.3) holds true (see, e.g., [6]). Throughout the paper, we use $u$ to denote this function.

The following result, which concerns the Malliavin derivatives of $(X, Y, Z)$ and provides another representation of $Z$, is due to Pardoux and Peng [9] (or see [6]).

LEMMA 2.1. *Assume* (A1), (A3); *and that* $f \in C_b^{0,1}, g \in C_b^1(\mathbb{R})$. *Let* $(\nabla X, \nabla Y, \nabla Z)$ *denote the solution to the following linear SDEs:*

$$(2.4) \quad \begin{aligned} \nabla X_t &= 1 + \int_0^t b_x(r, X_r)\nabla X_r\,dr + \int_0^t \sigma_x(r, X_r)\nabla X_r\,dW_r; \\ \nabla Y_t &= g'(X_T)\nabla X_T \\ &\quad + \int_t^T [f_x \nabla X_r + f_y \nabla Y_r + f_z \nabla Z_r]\,dr - \int_t^T \nabla Z_r\,dW_r. \end{aligned}$$

*Then it holds that, for $t \leq r \leq T$,*

$$
\begin{aligned}
D_t X_r &= \nabla X_r [\nabla X_t]^{-1} \sigma(t, X_t); \\
D_t Y_r &= \nabla Y_r [\nabla X_t]^{-1} \sigma(t, X_t); \\
D_t Z_r &= \nabla Z_r [\nabla X_t]^{-1} \sigma(t, X_t),
\end{aligned}
\tag{2.5}
$$

*where $D$ is the Malliavin derivative operator. Moreover, $u \in C_b^{0,1}$ and it holds that*

$$
u_x(t,x) = \nabla Y_t^{t,x}; \qquad Z_t = D_t Y_t = \nabla Y_t [\nabla X_t]^{-1} \sigma(t, X_t).
\tag{2.6}
$$

We note that Lemma 2.1 relies heavily on the differentiability of $g$. The next lemma, which gives a Feynman–Kac type representation formula of $u_x$, assumes instead that $\sigma$ is nondegenerate. The proof can be found in [6].

LEMMA 2.2. *Assume (A1), (A3), $\sigma \geq \frac{1}{K}$ and that $g$ is Lipschitz continuous. Then $u \in C_b^{0,1}([0,T) \times \mathbb{R})$, and (1.3) holds true in $[0,T) \times \mathbb{R}$. Moreover, we have the following representation formula of $u_x$ (and thus of $Z$):*

$$
u_x(t,x) = E^{t,x}\left\{ g(X_T)\bar{N}_T^t + \int_t^T f(r, X_r, Y_r, Z_r)\bar{N}_r^t \, dr \right\},
\tag{2.7}
$$

*where the superscript $^{t,x}$ indicates that the processes $(X,Y,Z)$ under expectation are solutions to $(2.3)$ [instead of $(1.1)$], and*

$$
\bar{N}_r^t \triangleq \frac{1}{r-t} \int_t^r \sigma^{-1}(s, X_s) \nabla X_s \, dW_s \, [\nabla X_t]^{-1}.
\tag{2.8}
$$

The following estimates are easy to prove (see, e.g., [2]).

LEMMA 2.3. *Assume that $\widetilde{b}, \widetilde{\sigma}: \Omega \times [0,T] \times \mathbb{R} \mapsto \mathbb{R}$ and $\widetilde{f}: \Omega \times [0,T] \times \mathbb{R}^2 \mapsto \mathbb{R}$ are $\mathbf{F}$-adapted random fields, such that they are uniformly Lipschitz continuous with respect to the spatial variable(s) and*

$$
E\left\{ \int_0^T [|\widetilde{b}(t,0)|^2 + |\widetilde{\sigma}(t,0)|^2 + |\widetilde{f}(t,0,0)|^2] \, dt \right\} < \infty.
$$

*For any $\xi \in L^2(\mathcal{F}_T)$, denote $(X,Y,Z)$ to be the solution to the following SDEs:*

$$
\begin{aligned}
X_t &= x + \int_0^t \widetilde{b}(s, X_s) \, ds + \int_0^t \widetilde{\sigma}(s, X_s) \, dW_s; \\
Y_t &= \xi + \int_t^T \widetilde{f}(s, Y_s, Z_s) \, ds - \int_t^T Z_s \, dW_s.
\end{aligned}
$$

*Then, for any $p \geq 2$, there exists a constant $C_p > 0$, depending only on $T, p$ and the Lipschitz constants of $\widetilde{b}$, $\widetilde{\sigma}$, $\widetilde{f}$, such that*

$$(2.9) \quad \begin{aligned} E\Big\{\sup_{0\leq t\leq T} |X_t|^p\Big\} \\ \leq C_p E\Big\{|x|^p + \int_0^T [|\widetilde{b}(t,0)|^p + |\widetilde{\sigma}(t,0)|^p]\,dt\Big\}; \end{aligned}$$

$$(2.10) \quad \begin{aligned} E\Big\{\sup_{0\leq t\leq T} |Y_t|^p + \Big(\int_0^T |Z_t|^2\,dt\Big)^{p/2}\Big\} \\ \leq C_p E\Big\{|\xi|^p + \int_0^T |\widetilde{f}(t,0,0)|^p\,dt\Big\}. \end{aligned}$$

We end this section with the exponential inequality (see, e.g., [10]).

LEMMA 2.4. *Assume $M$ is a continuous local martingale vanishing at $0$. Let $M_T^* \triangleq \sup_{0\leq t\leq T} |M_t|$ and $[M]_T$ denote its quadratic variation. Then for any $x, y > 0$,*

$$P(M_T^* > x, [M]_T < y) \leq \exp\Big(\frac{-x^2}{2y}\Big).$$

**3. Two "good" cases.** In this section we study two cases which generalize Lemmas 2.1 and 2.2, respectively. The first one assumes that $g$ is differentiable. Since the proof is more or less standard, we will just sketch it.

THEOREM 3.1. *Assume (A1), (A3), (A4) and that $f \in C_b^{0,1}, g \in C^1(\mathbb{R})$ such that $|g'(x)| \leq C\psi(x)$, where $\psi$ is defined in (A4). Then, for $\forall (t,x) \in [0,T] \times \mathbb{R}$:*

  (i) *$u \in C^{0,1}$, and (1.3) holds true;*
  (ii) *the following representation holds true:*

$$(3.1) \quad u_x(t,x) = E^{t,x}\Big\{g'(X_T)\nabla X_T + \int_t^T [f_x \nabla X_r + f_y \nabla Y_r + f_z \nabla Z_r]\,dr\Big\};$$

  (iii) *$|u_x(t,x)| \leq C\psi(x)$.*

PROOF. First, if $|g'(x)| \leq C$, then one gets (i) and (ii) immediately from Lemma 2.1. In general, by standard approximating arguments, one can prove (i) and (ii). Finally, by (3.1), (2.10) and (2.1), one can easily show that $|u_x(t,x)| = |\nabla Y_t^{t,x}| \leq C\psi(x)$. □

The next result is an extension of Lemma 2.2. We do not require $g$ to be continuous. Instead, we assume that $\sigma(T, \cdot)$ is nondegenerate.

THEOREM 3.2. *Assume* (A1), (A3), (A4) *and that there exists* $\delta_0 > 0$ *such that* $\frac{1}{K} \leq |\sigma(t,x)| \leq K$ *for* $\forall (t,x) \in [T - \delta_0, T] \times \mathbb{R}$. *Then:*

(i) $u \in C^{0,1}([0,T) \times \mathbb{R})$, *and* (1.3) *holds true;*
(ii) *the following estimates hold true: for* $0 \leq t \leq T$,

$$|u(t,x)| \leq C\psi(x); \qquad |u_x(t,x)| \leq \frac{C\psi(x)}{\sqrt{T-t}}.$$

PROOF. First, by Lemma 2.3, one can easily prove $|u(t,x)| \leq C\psi(x)$.

We now prove (i) and estimate $u_x$. We proceed in three steps.

*Step* 1. We restrict $t \in [T - \delta_0, T]$ and assume that $g \in C_b^1$. Then by Lemma 2.2, obviously, (i) holds true in $[T - \delta_0, T) \times \mathbb{R}$. Moreover, the representation formula (2.7) holds true. That is,

$$(3.2) \quad u_x(t,x) = E^{t,x}\bigg\{g(X_T)\bar{N}_T^t + \int_t^T f(r, X_r, u, u_x\sigma)(r, X_r)\bar{N}_r^t\,dr\bigg\}.$$

We shall use (3.2) to estimate $u_x$. We note that in the following the constant $C$ will *not* depend on the upper bound of $g'$. Define

$$A_t \stackrel{\triangle}{=} \sqrt{T-t}\sup_x \frac{|u_x(t,x)|}{\psi(x)}; \qquad B_t \stackrel{\triangle}{=} \sup_{t \leq s \leq T} A_s.$$

Recalling (2.8), one can check directly that, for $T - \delta_0 \leq t < r \leq T$, $E\{|\bar{N}_r^t|^p\} \leq C(r-t)^{-p/2}$. Note that $1 + |x| \leq C\psi(x)$, $|u(t,x)| \leq C\psi(x)$, and $|\sigma(t,x)| \leq K$. Then by (3.2), we have

$$|u_x(t,x)| \leq CE^{t,x}\bigg\{\psi(X_T)|\bar{N}_T^t|$$

$$+ \int_t^T [1 + |X_r| + |u(r,X_r)| + |(u_x\sigma)(r,X_r)|]|\bar{N}_r^t|\,dr\bigg\}$$

$$\leq CE^{t,x}\bigg\{\psi(X_T)|\bar{N}_T^t| + \int_t^T \psi(X_r)\Big[1 + \frac{A_r}{\sqrt{T-r}}\Big]|\bar{N}_r^t|\,dr\bigg\}$$

$$\leq C\psi(x)\bigg(\frac{1}{\sqrt{T-t}} + \int_t^T \Big[1 + \frac{A_r}{\sqrt{T-r}}\Big]\frac{dr}{\sqrt{r-t}}\bigg)$$

$$\leq C\psi(x)\bigg(\frac{1}{\sqrt{T-t}} + B_t\int_t^T \frac{dr}{\sqrt{(T-r)(r-t)}}\bigg)$$

$$= C\psi(x)\bigg(\frac{1}{\sqrt{T-t}} + B_t\int_0^1 \frac{dr}{\sqrt{r(1-r)}}\bigg) \leq C_0\psi(x)\bigg(\frac{1}{\sqrt{T-t}} + B_t\bigg),$$

where the last equality is due to the substitution $r = t + (T-t)r'$. Thus, $A_t \leq C_0[1 + B_t\sqrt{T-t}]$. Obviously, $B_t$ is decreasing, so for $t \leq s \leq T$, $A_s \leq C_0[1 + B_s\sqrt{T-s}] \leq C_0[1 + B_t\sqrt{T-t}]$. Therefore, $B_t \leq C_0[1 + B_t\sqrt{T-t}]$.

Without loss of generality, we assume $\delta_0 < (2C_0)^{-2}$. Then $C_0\sqrt{T-t} \leq \frac{1}{2}$, and, thus, $B_t \leq 2C_0$. This obviously implies that $|u_x(t,x)| \leq \frac{2C_0\psi(x)}{\sqrt{T-t}}$. We note again that $C_0$ does not depend on the bound of $g'$.

*Step* 2. We now assume that $g$ satisfies only (A4), but still restrict $t \in [T-\delta_0, T]$. One can easily find $g_n \in C_b^1$ such that $|g_n(x)| \leq 1 + \psi(x)$ for the same function $\psi$ and $\lim_{n\to\infty} g_n(x) = g(x)$ for $dx$-a.s. $x \in \mathbb{R}$, where $dx$ denotes the Lebesgue measure on $\mathbb{R}$. Let $(Y^n, Z^n)$ denote the solution to the BSDE:

$$Y_t^n = g_n(X_T) + \int_t^T f(r, X_r, Y_r^n, Z_r^n)\, dr - \int_t^T Z_r^n\, dW_r,$$

and then define $u^n(t,x) \triangleq Y_t^{n,t,x}$. Since $g_n \in C_b^1$, by the above arguments, we have $Z_t^n = (u_x^n \sigma)(t, X_t)$, the representation formula (2.7) holds true for $u_x^n$, and $|u_x^n(t,x)| \leq \frac{C\psi(x)}{\sqrt{T-t}}$, where $C$ is independent of $n$.

Noting that $\sigma(t,x) \geq \frac{1}{K}$, $X_T$ is absolutely continuous with respect to $dx$ (see, e.g., [8]). Thus, $\lim_{n\to\infty} g_n(X_T) = g(X_T)$, $P$-a.s. Then by standard arguments (see, e.g., [2]), one can show that

$$(3.3) \quad \lim_{n\to\infty} E\left\{\sup_{T-\delta_0 \leq t \leq T} |Y_t^n - Y_t|^2 + \int_{T-\delta_0}^T |Z_t^n - Z_t|^2\, dt\right\} = 0,$$

which implies that $\lim_{n\to\infty} u^n(t,x) = u(t,x)$. Moreover, recalling that

$$|Y_t| \leq C\psi(X_t), \qquad |Z_t| \leq \frac{C}{\sqrt{T-t}} \psi(X_t)|\sigma(t, X_t)|,$$

applying the dominated convergence theorem, one gets that $\lim_{n\to\infty} u_x^n(t,x) = v(t,x)$, where

$$v(t,x) = E^{t,x}\left\{g(X_T)\bar{N}_T^t + \int_t^T f(r, X_r, Y_r, Z_r)\bar{N}_r^t\, dr\right\}.$$

Obviously, $|v(t,x)| \leq \frac{C\psi(x)}{\sqrt{T-t}}$ and $Z_t = (v\sigma)(t, X_t)$. It remains to show that $v$ is continuous, and $u_x = v$. To this end, we note that, for any $\varepsilon > 0$, there exists an open set $O_\varepsilon \subset \mathbb{R}$ and a continuous function $g_\varepsilon$ such that: (i) the Lebesgue measure $|O_\varepsilon| \leq \varepsilon$; (ii) $g_\varepsilon(x) = g(x)$ for all $x \notin O_\varepsilon$; and (iii) $|g_\varepsilon(x)| \leq 1 + \psi(x)$. Denote

$$(3.4) \quad v_\varepsilon(t,x) = E^{t,x}\left\{g_\varepsilon(X_T)\bar{N}_T^t + \int_t^T f(r, X_r, Y_r, Z_r)\bar{N}_r^t\, dr\right\}.$$

We note that in (3.4) $(Y, Z)$ is still the solution to the BSDE with terminal value $g(X_T)$ [not $g_\varepsilon(X_T)$!]. Then

$$|v_\varepsilon(t,x) - v(t,x)| = |E^{t,x}\{[g_\varepsilon(X_T) - g(X_T)]\bar{N}_T^t\}|$$
$$\leq E^{t,x}\{[|g_\varepsilon(X_T)| + |g(X_T)|]|\bar{N}_T^t|; X_T \in O_\varepsilon\}$$
$$\leq \frac{C\psi(x)}{\sqrt{T-t}}\sqrt{P(X_T^{t,x} \in O_\varepsilon)}.$$

Again, by [8], one can easily show that $X_T^{t,x}$ has a bounded density function with respect to $dx$ (the bound may depend on $\frac{1}{T-t}$, however). Thus,

$$|v_\varepsilon(t,x) - v(t,x)| \leq C(T-t,x)\sqrt{|O_\varepsilon|} \leq C(T-t,x)\sqrt{\varepsilon}.$$

Here $C(T-t,x)$ is a constant depending on $T-t$ and $x$. Now for any $(t,x) \in [T-\delta_0, T] \times \mathbb{R}$ and $(t_n, x_n) \in [T-\delta_0, T] \times \mathbb{R}$ such that $\lim_{n\to\infty}(t_n, x_n) = (t,x)$, we have

$$|v(t_n, x_n) - v(t,x)|$$
$$\leq |v(t_n, x_n) - v_\varepsilon(t_n, x_n)| + |v_\varepsilon(t_n, x_n) - v_\varepsilon(t,x)| + |v_\varepsilon(t,x) - v(t,x)|$$
$$\leq [C(T-t_n, x_n) + C(T-t,x)]\sqrt{\varepsilon} + |v_\varepsilon(t_n, x_n) - v_\varepsilon(t,x)|.$$

Since $g_\varepsilon$ is continuous, by standard arguments, one can show that $\lim_{n\to\infty} v_\varepsilon(t_n, x_n) = v_\varepsilon(t,x)$. Thus, $\limsup_{n\to\infty} |v(t_n, x_n) - v(t,x)| \leq C\sqrt{\varepsilon}$. Since $\varepsilon$ is arbitrary, we have

$$\lim_{n\to\infty} |v(t_n, x_n) - v(t,x)| = 0.$$

That is, $v$ is continuous.

Finally, for any $(t,x) \in [T-\delta_0, T) \times \mathbb{R}$, note that

$$u^n(t,x) = u^n(t,0) + \int_0^x u_x^n(t,y)\,dy.$$

Let $n$ tend to $\infty$ and apply the dominated convergence theorem, we have

$$u(t,x) = u(t,0) + \int_0^x v(t,y)\,dy.$$

Since $v$ is continuous, we know that $u \in C^{0,1}$ and $u_x = v$. Therefore, for $t \in [T-\delta_0, T)$, (i) holds true and $|u_x(t,x)| \leq \frac{C\psi(x)}{\sqrt{T-t}}$.

*Step* 3. For $t \in [0, T-\delta_0]$, one may consider (1.1) as an FBSDE over $[0, T-\delta_0]$ such that the BSDE has terminal value $u(T-\delta_0, X_{T-\delta_0})$. By step 2, we have $|u_x(T-\delta_0, x)| \leq \frac{C}{\sqrt{\delta_0}}\psi(x)$. Then applying Theorem 3.1, we know (i) holds true in $[0, T-\delta_0]$ and $|u_x(t,x)| \leq \frac{C}{\sqrt{\delta_0}}\psi(x) \leq \frac{C\psi(x)}{\sqrt{T-t}}$. □

REMARK 3.3. All the results in this section hold true for high-dimensional FBSDEs.

**4. The case $f = 0$.** In this section we study the case that $f = 0$. In this case the BSDE in (1.1) becomes

(4.1) $$Y_t = g(X_T) - \int_t^T Z_r\,dW_r.$$

4.1. *A counterexample.* Note that both in Theorem 3.1 and in Theorem 3.2, $u_x(t,x)$ exists for all $t < T$. When $\sigma$ is degenerate and $g$ is not Lipschitz continuous, however, $u_x(t,x)$ may not exist for some $t < T$. An obvious example is that $\sigma = b = 0$ and $g = \mathbb{1}_{\{x>0\}}$, then $u(t,x) = \mathbb{1}_{\{x>0\}}$ for $\forall t \in [0,T]$, and, hence, $u$ is discontinuous in $x$. But, since $\sigma = 0$, in light of (1.3) in this example, one may still view $Z_t = 0$. In fact, $Y_t = \mathbb{1}_{\{X_t > 0\}}$, $Z_t = 0$ are indeed the solution to (4.1).

The following example shows that $Z_t$ may also blow up for some $t < T$.

EXAMPLE 1. *Assume $\alpha \in (0,1)$ and $\beta \in (0, \frac{\alpha}{2(1-\alpha)})$. Let $T = 2$, and*

$$g(x) \triangleq \frac{x}{|x|^\alpha}; \qquad b(t,x) \triangleq 0 \qquad \forall t \in [0,2];$$

$$\sigma(t,x) \triangleq \begin{cases} (1-t)^\beta, & t \in [0,1], \\ 0, & t \in (1,2]. \end{cases}$$

*Then $Z_t$ may blow up at $t = 1$.*

PROOF. For $t \in [1,2]$, obviously $u(t,x) = g(x)$ and $Z_t = 0$. For $t \in [0,1)$, we have

$$u(t,x) = E^{t,x}\{u(1,X_1)\} = E\left\{g\left(x + \int_t^1 (1-s)^\beta \, dW_s\right)\right\}.$$

Note that $\int_t^1 (1-s)^\beta \, dW_s$ has normal distribution with mean 0 and variance $\sigma_0^2 \triangleq \int_t^1 (1-s)^{2\beta} \, ds = \frac{1}{1+2\beta}(1-t)^{1+2\beta}$. Thus,

$$u(t,x) = \frac{1}{\sqrt{2\pi}\sigma_0} \int_\mathbb{R} g(y) \exp\left(-\frac{(y-x)^2}{2\sigma_0^2}\right) dy.$$

Therefore,

$$u_x(t,0) = \frac{1}{\sqrt{2\pi}\sigma_0} \int_\mathbb{R} g(y) \frac{y}{\sigma_0^2} e^{-y^2/(2\sigma_0^2)} \, dy = \frac{1}{\sqrt{2\pi}\sigma_0^3} \int_\mathbb{R} |y|^{2-\alpha} e^{-y^2/(2\sigma_0^2)} \, dy.$$

By using the substitution $y = \sigma_0 y'$, we get

$$u_x(t,0) = C\sigma_0^{-\alpha} = C(1-t)^{-\alpha(1+2\beta)/2}.$$

Thus,

$$(u_x \sigma)(t,0) = C(1-t)^{-\alpha(1+2\beta)/2 + \beta} = C(1-t)^{\beta(1-\alpha) - \alpha/2}.$$

Since $\beta \in (0, \frac{\alpha}{2(1-\alpha)})$, we have $\lim_{t \uparrow 1}(u_x\sigma)(t,0) = \infty$. That implies that, if $X_1 = 0$, then $Z_t \to \infty$ as $t \uparrow 1$. □

4.2. *Location of discontinuous points.* Note that in Example 1, $Z_t$ blows up only at $t = 1$. In fact, for quite general FBSDEs, along each path $Z_t$ is discontinuous at most at one point. In this section we locate this possible discontinuous point and we shall prove later that $Z_t$ is continuous elsewhere.

To this end, we introduce the following notation. For $\forall (t, x) \in [0, T) \times \mathbb{R}$, let $\eta_{\cdot}^{t,x}$ be the (deterministic) characteristic of the FSDE in (2.3). That is, $\eta_{\cdot}^{t,x}$ is the solution to the following integral equation:

$$(4.2) \qquad \eta_s^{t,x} = x + \int_t^s b(r, \eta_r^{t,x}) \, dr.$$

As usual, we omit the superscript $^{0,x}$ when $t = 0$. Define

$$(4.3) \qquad \begin{aligned} \Gamma^0 &\triangleq \left\{ (t, x) : \max_{t \leq s \leq T} |\sigma(s, \eta_s^{t,x})| > 0 \right\}; \\ \Gamma^n &\triangleq \left\{ (t, x) : \max_{t \leq s \leq T} |\sigma(s, \eta_s^{t,x})| \geq \frac{1}{n} \right\}; \end{aligned}$$

$$(4.4) \qquad \begin{aligned} \Gamma &\triangleq \{(t, x) : |\sigma(t, x)| > 0\}; \\ \tau &\triangleq \inf\{0 \leq t \leq T : (t, X_t) \notin \Gamma^0\}. \end{aligned}$$

Note that $\max_{t \leq s \leq T} |\sigma(s, \eta_s^{t,x})|$ is continuous in $(t, x)$, then the following results are obvious.

LEMMA 4.1. *Assume* (A1). *Then:*

(i) $\Gamma \subset \Gamma^0 = \bigcup_{n=1}^\infty \Gamma^n$, and $\Gamma, \Gamma^n$ are open for $\forall n \geq 0$;
(ii) For $\forall (t, x) \in \Gamma^0$, $\sup_{s \in [t,T]} |\sigma(s, X_s^{t,x})| > 0$, a.s.;
(iii) $\tau$ is a stopping time;
(iv) $(t, X_t) \in \Gamma^0$ for all $t < \tau$;
(v) $(t, X_t) \notin \Gamma^0$ and $\sigma(t, X_t) = 0$ for all $t \geq \tau$.

4.3. *Representation formula.* In this section we formally derive a new Feynman–Kac type representation formula for $u_x(t, x)$ (and, hence, for $Z_t$) in $\Gamma^0$. We shall follow the arguments in [6].

Fix $(t, x) \in \Gamma^0$. We first assume that (A1), (A2) hold true and that $g \in C_b^1(\mathbb{R})$. Note that $u(t, x) = E^{t,x}\{g(X_T)\}$ and that $\nabla X_t$ is the derivative of the flow $X_t^x$ with respect to initial value $x$ (see [9]). Then

$$(4.5) \qquad u_x(t, x) = E^{t,x}\{g'(X_T) \nabla X_T\}.$$

For $t < s < T$, by Lemma 2.1, we have

$$D_s g(X_T) = g'(X_T) D_s X_T = g'(X_T) \nabla X_T [\nabla X_s]^{-1} \sigma(s, X_s).$$

Multiply both sides by $\nabla X_s \sigma(s, X_s)$ and integrate over $[t, T]$, we get

$$\int_t^T \nabla X_s \sigma(s, X_s) D_s g(X_T) \, ds = g'(X_T) \nabla X_T \int_t^T \sigma^2(s, X_s) \, ds.$$

Denote

(4.6) $$\gamma_t \triangleq \sigma(t, X_t); \qquad \Lambda_r^t \triangleq \int_t^r \gamma_s^2 \, ds.$$

Since $(t, x) \in \Gamma^0$, by Lemma 4.1, we have $\Lambda_T^t > 0$, a.s. Then

$$g'(X_T) \nabla X_T = \int_t^T \frac{\gamma_s \nabla X_s}{\Lambda_T^t} D_s g(X_T) \, ds.$$

Recalling (4.5) and applying the integration by parts formula of Malliavin calculus, we have

$$u_x(t, x) = E^{t,x} \bigg\{ g(X_T) \int_t^T \frac{1}{\Lambda_T^t} \gamma_s \nabla X_s \, dW_s \bigg\}$$

$$= E^{t,x} \bigg\{ g(X_T) \bigg[ \frac{1}{\Lambda_T^t} \int_t^T \gamma_s \nabla X_s \, dW_s - \int_t^T D_s \bigg( \frac{1}{\Lambda_T^t} \bigg) \gamma_s \nabla X_s \, ds \bigg] \bigg\}.$$

Note that

$$D_s \bigg( \frac{1}{\Lambda_T^t} \bigg) \gamma_s \nabla X_s = -\frac{D_s \Lambda_T^t}{|\Lambda_T^t|^2} \gamma_s \nabla X_s$$

$$= -\frac{2}{|\Lambda_T^t|^2} \int_s^T \sigma_x(r, X_r) \gamma_r D_s X_r \, dr \, \gamma_s \nabla X_s$$

$$= -\frac{2}{|\Lambda_T^t|^2} \int_s^T \sigma_x(r, X_r) \gamma_r \nabla X_r \, dr \, \gamma_s^2.$$

Thus,

$$u_x(t, x) = E^{t,x} \bigg\{ g(X_T) \bigg[ \frac{1}{\Lambda_T^t} \int_t^T \gamma_s \nabla X_s \, dW_s$$

$$+ \frac{2}{|\Lambda_T^t|^2} \int_t^T \int_s^T \sigma_x(r, X_r) \gamma_r \nabla X_r \, dr \, \gamma_s^2 \, ds \bigg] \bigg\}$$

$$= E^{t,x} \bigg\{ g(X_T) \bigg[ \frac{1}{\Lambda_T^t} \int_t^T \gamma_s \nabla X_s \, dW_s$$

$$+ \frac{2}{|\Lambda_T^t|^2} \int_t^T \sigma_x(r, X_r) \gamma_r \nabla X_r \int_t^r \gamma_s^2 \, ds \, dr \bigg] \bigg\}.$$

Denote

(4.7) $$N_r^t \triangleq \frac{1}{\Lambda_r^t} \bigg[ \int_t^r \gamma_s \nabla X_s \, dW_s + 2 \int_t^r \frac{\Lambda_s^t}{\Lambda_r^t} \sigma_x(s, X_s) \gamma_s \nabla X_s \, ds \bigg] [\nabla X_t]^{-1}.$$

Noting that $\nabla X_t^{t,x} = 1$, we have

(4.8) $$u_x(t,x) = E^{t,x}\{g(X_T)N_T^t\}.$$

REMARK 4.2. In (2.8) one has to assume $\sigma(t,x) \neq 0$ for all $(t,x)$. But in (4.7), even if $\sigma(t,x) = 0$, one can still define $N_r^{t,x}$ as long as $\Lambda_r^t > 0$.

4.4. *Estimate of* $[\Lambda_r^t]^{-1}$. In light of (4.7) and (4.8), obviously one needs to estimate $[\Lambda_r^t]^{-1}$. We have the following result.

LEMMA 4.3. *Assume* (A1), (A2) *and* (A5). *If* $(t,x) \in \Gamma^n$ *for some* $n \geq 1$, *then for any* $p \geq 1$, *there exists a constant* $C_p$ *depending on* $p, K$ *and* $T$, *such that*

$$E^{t,x}\{|\Lambda_r^t|^{-p}\} \leq C_p[n^{2p}\delta_n^{-p} + n^{3p+1}\delta_n],$$

*where* $\delta_n \triangleq \delta(\frac{1}{4n}) \wedge \frac{r-t}{2}$ *and* $\delta(\cdot)$ *is defined in* (2.2).

PROOF. Without loss of generality, we prove the result only for $t = 0$ and $r = T$. To this end, we fix $(0,x) \in \Gamma^n, n \geq 1$, and omit the superscripts $^{0,x}$ whenever there is no confusion. So it suffices to prove that

(4.9) $$E\{\Lambda_T^{-p}\} \leq C_p[n^{2p}\delta_n^{-p} + n^{3p+1}\delta_n],$$

where $\delta_n \triangleq \delta(\frac{1}{4n}) \wedge \frac{T}{2}$. Note that

(4.10) $$E\{\Lambda_T^{-p}\} = p\int_0^\infty P\{\Lambda_T < u\}u^{-p-1}\,du.$$

We shall estimate $P\{\Lambda_T < u\}$ for small $u$ below.

Since $(0,x) \in \Gamma^n$, there exists $T_0 \in [0,T]$ such that $|\sigma(T_0, \eta_{T_0})| > \frac{1}{n}$. Note that

(4.11) $$P\{\Lambda_T < u\} = I_1(u) + I_2(u),$$

where

$$I_1(u) \triangleq P\left\{\Lambda_T < u, |\gamma_{T_0}| > \frac{1}{2n}\right\}; \qquad I_2(u) \triangleq P\left\{\Lambda_T < u, |\gamma_{T_0}| \leq \frac{1}{2n}\right\}.$$

We first estimate $I_2(u)$. Denote

(4.12) $$\Delta X_t \triangleq X_t - \eta_t.$$

Then one can easily get

(4.13) $$I_2(u) \leq P\left\{\Lambda_T < u, |\Delta X_{T_0}| \geq \frac{1}{2Kn}\right\} \leq P\left\{\Lambda_{T_0} < u, |\Delta X_{T_0}| \geq \frac{c}{n}\right\}.$$

Note that
$$\Delta X_t = \int_0^t \gamma_s \, dW_s + \int_0^t \beta_s \Delta X_s \, ds,$$
where $\beta_s \triangleq \frac{b(s,X_s)-b(s,\eta_s)}{\Delta X_s}$ is bounded. Denote

(4.14) $$L_t \triangleq \exp\left(-\int_0^t \beta_s \, ds\right); \qquad M_t \triangleq \int_0^t L_s \gamma_s \, dW_s.$$

Then $\Delta X_{T_0} = L_{T_0}^{-1} M_{T_0}$. Moreover, both $L_t$ and $L_t^{-1}$ are bounded, and $M_t$ is a martingale. Then one has
$$|\Delta X_{T_0}| \leq C|M_{T_0}|; \qquad [M]_{T_0} = \int_0^{T_0} L_t^2 \gamma_t^2 \, dt \leq C\Lambda_{T_0}.$$

Now applying Lemma 2.4, we get from (4.13) that

(4.15) $$I_2(u) \leq P\left\{[M]_{T_0} < Cu, |M_{T_0}| \geq \frac{c}{n}\right\} \leq \exp\left(-\frac{c}{n^2 u}\right).$$

We next estimate $I_1(u)$. To this end, we recall that $\delta_n \triangleq \delta(\frac{1}{4n}) \wedge \frac{T}{2}$. Denote
$$\delta_u \triangleq 64n^2 u;$$
$$\Delta_u \triangleq [T_0 - \delta_u, T_0 + \delta_u] \cap [0,T];$$
$$\Delta_0 \triangleq [T_0 - \delta_n, T_0 + \delta_n] \cap [0,T].$$
For $u \leq \frac{\delta_n}{64n^2}$, obviously we have
$$\delta_u \leq \delta_n; \qquad \Delta_u \subset \Delta_0; \qquad |\Delta_u| \geq \delta_u.$$
Using the facts that
$$\Lambda_T \geq \int_{\Delta_u} \gamma_t^2 \, dt \geq |\Delta_u| \inf_{t \in \Delta_u} |\gamma_t|^2 \geq \delta_u \inf_{t \in \Delta_u} |\gamma_t|^2 = 64n^2 u \inf_{t \in \Delta_u} |\gamma_t|^2,$$
and that
$$|\sigma(t, X_{T_0}) - \sigma(T_0, X_{T_0})| \leq \frac{1}{4n} \qquad \forall t \in \Delta_u \subset \Delta_0,$$
we have
$$I_1(u) = P\left\{\Lambda_T < u, |\gamma_{T_0}| > \frac{1}{2n}, \inf_{t \in \Delta_u} |\gamma_t| < \frac{1}{8n}\right\}$$
$$\leq P\left\{\Lambda_T < u, \sup_{t \in \Delta_u} |\sigma(t, X_t) - \sigma(t, X_{T_0})| \geq \frac{1}{8n}\right\}$$
$$\leq P\left\{\Lambda_T < u, \sup_{t \in \Delta_u} |X_t - X_{T_0}| \geq \frac{1}{8Kn}\right\}.$$

Since
$$X_{T_0} - X_t = \int_t^{T_0} b(s, X_s)\, ds + \int_t^{T_0} \gamma_s\, dW_s,$$

we have

(4.16) $$I_1(u) \leq I_{11}(u) + I_{12}(u),$$

where

$$I_{11}(u) \triangleq P\bigg\{\int_{\Delta_u} |b(t, X_t)|\, dt \geq \frac{1}{16Kn}\bigg\};$$
$$I_{12}(u) \triangleq P\bigg\{\Lambda_T < u,\, \sup_{t \in \Delta_u}\bigg|\int_t^{T_0} \gamma_s\, dW_s\bigg| \geq \frac{1}{16Kn}\bigg\}.$$

Note that

(4.17)
$$I_{11}(u) \leq C_p n^{p+1} E\bigg\{\bigg(\int_{\Delta_u} |b(t, X_t)|\, dt\bigg)^{p+1}\bigg\}$$
$$\leq C_p (n\delta_u)^{p+1} = C_p n^{3p+3} u^{p+1}.$$

Moreover, analogous to (4.15), by applying Lemma 2.4, one can show that

(4.18) $$I_{12}(u) \leq \exp\bigg(-\frac{c}{n^2 u}\bigg).$$

Plugging (4.17) and (4.18) into (4.16), we get

$$I_1(u) \leq C_p n^{3p+3} u^{p+1} + \exp\bigg(-\frac{c}{n^2 u}\bigg).$$

Noting that $\exp(-\frac{c}{n^2 u}) u^{-p-1}$ takes its maximum value at $u = \frac{c}{n^2(p+1)}$, the above inequality, together with (4.15), (4.11) and (4.10), implies that

$$E\{\Lambda_T^{-p}\} \leq C_p \bigg[\int_0^{\delta_n/(64n^2)} [n^{3p+3} u^{p+1} + e^{-c/(n^2 u)}] u^{-p-1}\, du$$
$$+ \int_{\delta_n/(64n^2)}^{\infty} u^{-p-1}\, du\bigg]$$
$$\leq C_p\bigg[n^{3p+1}\delta_n + n^{2(p+1)}\frac{\delta_n}{n^2} + (n^2 \delta_n^{-1})^p\bigg] \leq C_p[n^{3p+1}\delta_n + (n^2 \delta_n^{-1})^p],$$

which proves (4.9), and whence the lemma. □

4.5. *Main results.* Now we are ready to state the main result of this section.

THEOREM 4.4. *Assume* (A1), (A2), (A4), (A5) *and* $f = 0$. *Then:*

(i) $u \in C^{0,1}(\Gamma^0)$;
(ii) $Z_t = \begin{cases} (u_x\sigma)(t, X_t), & t < \tau; \\ 0, & t \geq \tau; \end{cases}$
(iii) $u_x(t, x) = E^{t,x}\{g(X_T)N_T^t\} \quad \forall (t, x) \in \Gamma^0$;
(iv) $|u_x(t, x)| \leq \dfrac{C_n \psi(x)}{\sqrt{T-t}} \quad \forall (t, x) \in \Gamma^n$.

REMARK 4.5. By Lemma 4.1(iv), we have $(t, X_t) \in \Gamma^0$ for $t < \tau$. Thus, by (i) and (ii) of Theorem 4.4, $Z_t$ is continuous along each path except possibly at $t = \tau$. But as we see in Example 1, $Z_t$ may have no finite left limit at $t = \tau$.

PROOF OF THEOREM 4.4. First, for $\forall (t, x) \in \Gamma^n$, by (4.7) and Lemma 4.3, one can easily prove that

$$E^{t,x}\{|N_T^t|^2\} \leq \frac{C_n}{T-t}. \tag{4.19}$$

Then (iv) is a direct consequence of (iii).

We now prove (i)–(iii). To this end, we first assume $g \in C_b^1(\mathbb{R})$. Then by applying Theorem 3.1, we know that $u \in C^{0,1}$ and $Z_t = (u_x\sigma)(t, X_t)$ for all $t \in [0, T]$. Especially, for $t \geq \tau$, by Lemma 4.1(v), we know $\sigma(t, X_t) = 0$, and, thus, $Z_t = 0$. Moreover, for $\forall (t, x) \in \Gamma^0$, by Lemma 4.1(i), there exists some $n \geq 1$ such that $(t, x) \in \Gamma^n$. Then by (4.19), we know $E^{t,x}\{|g(X_T)N_T^t|\} < \infty$. Now (iii) follows the arguments in Section 4.3.

In general case, that is, $g$ satisfies only (A4), we follow the arguments in step 2 of the proof for Theorem 3.2. Let $g_m \in C_b^1(\mathbb{R})$ such that $|g_m(x)| \leq \psi(x)$ and $\lim_{m\to\infty} g_m(x) = g(x)$ for $dx$-a.s. $x \in \mathbb{R}$. Define $(Y^m, Z^m)$ and $u^m$ as in Theorem 3.2. Note that $\Gamma^0, \Gamma^n$ and $\tau$ are independent of $g$, and so is $N$. Now for $\forall (t, x) \in \Gamma^n$ and any $m$,

$$u_x^m(t, x) = E^{t,x}\{g_m(X_T)N_T^t\}; \qquad |u_x^m(t, x)| \leq \frac{C_n \psi(x)}{\sqrt{T-t}}.$$

By the dominated convergence theorem, we have

$$u_x^m(t, x) \to v(t, x) \stackrel{\triangle}{=} E^{t,x}\{g(X_T)N_T^t\}, \qquad m \to \infty.$$

By a line by line analogy of step 2 of the proof for Theorem 3.2, one can show that $v$ is continuous and $u_x(t, x) = v(t, x)$. That proves (i) and (iii). Moreover, $Z^m \to Z$, then (ii) holds true. □

**5. General case.**

5.1. *Main results.* In this section we investigate FBSDE (1.1) with nonlinear $f$. We shall modify some assumptions:

(A2′) $\sigma$ is uniformly Hölder-$\alpha$ continuous in $t$ for some $\alpha > \frac{1}{2}$.
(A3′) $f$ is independent of $z$, that is, $f(t,x,y,z) = f(t,x,y)$. Moreover, $f$ is continuous in $t$ and uniformly Lipschitz continuous in $x, y$.
(A3″) $f(t,x,y,z) = f_1(t,x,y) + f_2(t,x)z$, where $f_1, f_2$ are continuous and uniformly Lipschitz continuous in $x, y$. Moreover, $f_2$ is bounded.

The following result gives an important estimate for $u_x$.

THEOREM 5.1. *Assume* (A1), (A2′), (A3′), (A4), (A5) *and that* $f \in C_b^{0,1}$ *and* $g \in C_b^1$. *Then for any* $(t,x) \in \Gamma^n$, *we have*

$$|u_x(t,x)| \leq \frac{C_n \psi(x)}{\sqrt{T-t}},$$

*where* $C_n$ *depends on* $K, T, \alpha, \psi$ *and* $n$, *but does* not *depend on the upper bound of* $g'$.

The main result of the paper is the following theorem.

THEOREM 5.2. *Assume* (A1), (A2′), (A3″), (A4) *and* (A5). *Then:*

(i) $u \in C^{0,1}(\Gamma)$, *and for* $\forall (t,x) \in \Gamma$, *we have*

$$u_x(t,x) = E^{t,x}\left\{g(X_T)N_T^t + \int_t^T f(r, X_r, Y_r, Z_r) N_r^t \, dr\right\}.$$

(ii) $u$ *is locally Lipschitz continuous in* $x$ *in* $\Gamma^0$, *and there exists a constant* $C_n$ *depending on* $K, T, \alpha, \psi$ *and* $n$ *such that*

$$|u_x(t,x)| \leq \frac{C_n \psi(x)}{\sqrt{T-t}} \qquad \forall (t,x) \in \Gamma^n.$$

*Here* $u_x(t,x)$ *is understood as a generalized derivative if* $u$ *is not differentiable in* $x$ *at* $(t,x)$.

(iii) *Understanding* $u_x$ *as in* (ii), *we have*

$$Z_t = \begin{cases} (u_x \sigma)(t, X_t), & t < \tau, \\ 0, & t \geq \tau. \end{cases}$$

REMARK 5.3. $Z_t$ is continuous except possibly at $t = \tau$. In fact, for $t < \tau$, if $(t, X_t) \in \Gamma$, then $Z_t$ is continuous by Theorem 5.2(i). If $(t, X_t) \in \Gamma^0 \backslash \Gamma$, we have $Z_t = 0$ and by the estimate in Theorem 5.2(ii), we know $Z_t$ is also continuous.

REMARK 5.4. If we remove (A5) in Theorems 5.1 and 5.2, then analogously one can show that

$$|u_x(t,x)| \leq \frac{C_n}{\sqrt{T-t}}(1+|x|)\psi(x) \qquad \forall\, (t,x) \in \Gamma^n.$$

We do not know whether or not similar results will hold true if $f$ is nonlinear on $z$. The proof of Theorems 5.1 is quite lengthy. We split it into several lemmas.

5.2. *Fine estimates of* $[\Lambda_r^t]^{-1}$. We first prove two lemmas which improve Lemma 4.3.

LEMMA 5.5. *For any $p \geq 1$ and $0 < \mu < \frac{1}{2}$, there exists a constant $C_{p,\mu}$, depending only on $p, \mu$, such that for any $T$ and any square integrable process $\gamma_t$, it holds that*

$$E\left\{\left[\left(\sup_{0\leq t\leq T}|\gamma_t|\right)^{-1}\sup_{0\leq t_1<t_2\leq T}\frac{\int_{t_1}^{t_2}\gamma_t\,dW_t}{(t_2-t_1)^\mu}\right]^p\right\} \leq C_{p,\mu}T^{(1/2-\mu)p} < \infty,$$

*where the integrand $\frac{0}{0}$ is considered as 0.*

PROOF. We proceed in two steps.

*Step* 1. We prove the following estimate:

(5.1) $$E\left\{\left|\left(\sup_{0\leq t\leq T}|\gamma_t|\right)^{-1}\int_0^T \gamma_t\,dW_t\right|^p\right\} \leq C_p T^{p/2}.$$

We first assume $0 < c \leq |\gamma_t| \leq C < \infty$. Denote

$$M_t \triangleq \int_0^t \gamma_s\,dW_s; \qquad \gamma_t^* \triangleq \sup_{0\leq s\leq t}|\gamma_s|.$$

Note that $\gamma_t^*$ is increasing and $\gamma_t^2 \leq |\gamma_t^*|^2$. Then applying Itô's formula, we have

$$d\left(\frac{M_t}{\gamma_t^*}\right)^2 = \frac{2M_t\,dM_t + \gamma_t^2\,dt}{|\gamma_t^*|^2} - \frac{2M_t^2}{|\gamma_t^*|^3}d\gamma_t^* \leq \frac{2M_t}{|\gamma_t^*|^2}dM_t + dt.$$

Obviously, $\lim_{t\downarrow 0}\frac{M_t}{\gamma_t^*} = 0$, thus, $E\{(\frac{M_t}{\gamma_t^*})^2\} \leq t$. Similarly,

$$d\left(\frac{M_t}{\gamma_t^*}\right)^{2n} = n\left(\frac{M_t}{\gamma_t^*}\right)^{2(n-1)}d\left(\frac{M_t}{\gamma_t^*}\right)^2 + \frac{n(n-1)}{2}\left(\frac{M_t}{\gamma_t^*}\right)^{2(n-2)}\frac{4M_t^2}{|\gamma_t^*|^4}\gamma_t^2\,dt$$

$$\leq n\left(\frac{M_t}{\gamma_t^*}\right)^{2(n-1)}\frac{2M_t}{|\gamma_t^*|^2}dM_t + C_n\left(\frac{M_t}{\gamma_t^*}\right)^{2(n-1)}dt.$$

By induction, one can prove $E\{(\frac{M_t}{\gamma_t^*})^{2n}\} \leq C_n t^n$. Then (5.1) holds true for any $p \geq 1$.

When $\gamma$ is unbounded or degenerate, by standard truncation procedure, one can easily prove (5.1).

*Step* 2. We follow the proof of Kolmogorov's continuity criterion (see, e.g., [10], Theorem 1.2.1). First, similar to (5.1), one can easily show that, for any $0 \leq t_1 < t_2 \leq T$,

$$\begin{aligned}E\left\{\left|(\gamma_T^*)^{-1}\int_{t_1}^{t_2}\gamma_t\,dW_t\right|^p\right\} &\leq E\left\{\left|\left(\sup_{t_1\leq t\leq t_2}|\gamma_t|\right)^{-1}\int_{t_1}^{t_2}\gamma_t\,dW_t\right|^p\right\}\\ &\leq C_p(t_2-t_1)^{p/2}.\end{aligned} \tag{5.2}$$

For any integer $n \geq 1$, let $D_n \triangleq \{i2^{-n}T : i = 1, \ldots, 2^n\}$ and $D \triangleq \bigcup_n D_n$. Denote $K_n \triangleq \sup_{1 \leq i \leq 2^n} |M_{i2^{-n}T} - M_{(i-1)2^{-n}T}|$. For $p > 2$, by (5.2), we have

$$\begin{aligned}E\{|\gamma_T^*|^{-p}K_n^p\} &\leq \sum_{i=1}^{2^n} E\{|\gamma_T^*|^{-p}|M_{i2^{-n}T} - M_{(i-1)2^{-n}T}|^p\}\\ &\leq C\sum_{i=1}^{2^n} 2^{-np/2}T^{p/2} = C2^{-n(p-2)/2}T^{p/2}.\end{aligned} \tag{5.3}$$

Now for any $t_1, t_2 \in D$ such that $t_2 - t_1 \leq 2^{-n}T$, one can easily show that $|M_{t_2} - M_{t_1}| \leq 2\sum_{m=n}^{\infty} K_m$. Denote

$$M^* \triangleq \sup\left\{\frac{|M_{t_2}-M_{t_1}|}{(t_2-t_1)^\mu}; t_1, t_2 \in D\right\}.$$

Then

$$\begin{aligned}M^* &\leq \sum_{n=0}^{\infty} \sup_{2^{-(n+1)}T \leq |t_2-t_1| \leq 2^{-n}T}\left\{\frac{|M_{t_2}-M_{t_1}|}{(t_2-t_1)^\mu}; t_1,t_2\in D\right\}\\ &\leq \sum_{n=0}^{\infty} 2^{(n+1)\mu}T^{-\mu}\sup_{|t_2-t_1|\leq 2^{-n}T}\{|M_{t_2}-M_{t_1}|; t_1,t_2\in D\}\\ &\leq \sum_{n=0}^{\infty} 2^{1+(n+1)\mu}T^{-\mu}\sum_{m=n}^{\infty}K_m \leq CT^{-\mu}\sum_{n=1}^{\infty}2^{n\mu}K_n.\end{aligned}$$

Now for $p > \frac{2}{1-2\mu}$, we have

$$\begin{aligned}\||\gamma_T^*|^{-1}M^*\|_p &\leq CT^{-\mu}\sum_{n=1}^{\infty}2^{n\mu}\||\gamma_T^*|^{-1}K_n\|_p\\ &\leq CT^{-\mu}\sum_{n=1}^{\infty}2^{n\mu}2^{-n(p-2)/(2p)}T^{1/2}\end{aligned}$$

$$= C \sum_{n=1}^{\infty} 2^{-n(1/2-1/p-\mu)} T^{1/2-\mu} = CT^{1/2-\mu}.$$

Now applying Fatou's lemma, we prove the lemma for $p > \frac{2}{1-2\mu}$. Then by the Schwarz inequality, the lemma holds true for smaller $p$. □

From now on we use $\gamma_t$ to denote $\sigma(t, X_t)$ again.

LEMMA 5.6. *Assume* (A1), (A2) *and* (A5). *Then for any* $(t,x) \in \Gamma^0$ *and any* $\beta \geq 0$, *there exists a constant* $C$, *depending only on* $K$, *such that*

$$|\Lambda_r^{t,x}|^{-1} \leq C|\gamma_r^{t,x,*}|^{-2} \bigg[ \delta\bigg(\frac{(r-t)^\beta \gamma_r^{t,x,*}}{4}\bigg)^{-1}$$
$$+ (r-t)^{-\beta}|\gamma_r^{t,x,*}|^{-1} + \frac{|\xi_r^{t,x}|^3}{(r-t)^{-3\beta}} + \frac{\mathbb{1}_{G_r^{t,x}}}{r-t}\bigg],$$

*where*

(5.4)
$$\gamma_r^{t,x,*} \triangleq \sup_{t \leq s \leq r} |\gamma_s^{t,x}|; \qquad \xi_r^{t,x} \triangleq \frac{1}{\gamma_r^{t,x,*}} \sup_{t \leq s_1 < s_2 \leq r} \bigg|\frac{\int_{s_1}^{s_2} \gamma_s^{t,x} dW_s}{(s_2 - s_1)^{1/3}}\bigg|;$$
$$G_r^{t,x} \triangleq \bigg\{\sup_{t \leq s \leq r} |\gamma_s^{t,x} - \sigma(t,x)| \leq |\sigma(t,x)|(r-t)^\beta \bigg\}.$$

PROOF. Without loss of generality, we assume $t = 0, r = T$ and $T \leq 1$. In the following we omit the superscript $^{0,x}$. We shall prove that

(5.5) $$\Lambda_T^{-1} \leq C|\gamma_T^*|^{-2}\bigg[\delta\bigg(\frac{T^\beta \gamma_T^*}{4}\bigg)^{-1} + T^{-\beta}|\gamma_T^*|^{-1} + T^{-3\beta}\xi_T^3 + T^{-1}\mathbb{1}_{G_T}\bigg].$$

First, there exists $\tilde{\tau}$ such that $\gamma_{\tilde{\tau}} = \gamma_T^*$. We note that, in general, $\tilde{\tau}$ is *not* a stopping time. Denote

$$\tilde{\tau}_- \triangleq \sup\bigg\{t \in [0, \tilde{\tau}] : |\gamma_t| = \bigg(1 - \frac{T^\beta}{2}\bigg)\gamma_T^*\bigg\},$$
$$\tilde{\tau}_+ \triangleq \inf\bigg\{t \in [\tilde{\tau}, T] : |\gamma_t| = \bigg(1 - \frac{T^\beta}{2}\bigg)\gamma_T^*\bigg\},$$

where we take the convention that $\sup \phi = 0$ and $\inf \phi = T$ for the empty set $\phi$. Then for $\forall t \in [\tilde{\tau}_-, \tilde{\tau}_+]$, we have $|\gamma_t| \geq \frac{\gamma_T^*}{2}$. Therefore,

(5.6) $$\Lambda_T \geq \tfrac{1}{4}|\gamma_T^*|^2(\tilde{\tau}_+ - \tilde{\tau}_-).$$

So it suffices to estimate $(\tilde{\tau}_+ - \tilde{\tau}_-)^{-1}$.

To this end, we note that, if $\tilde{\tau}_- > 0$, then

$$\frac{T^\beta \gamma_T^*}{2} = |\gamma_{\tilde{\tau}} - \gamma_{\tilde{\tau}_-}| \leq |\sigma(\tilde{\tau}, X_{\tilde{\tau}}) - \sigma(\tilde{\tau}_-, X_{\tilde{\tau}})| + K|X_{\tilde{\tau}} - X_{\tilde{\tau}_-}|$$

$$\leq |\sigma(\tilde{\tau}, X_{\tilde{\tau}}) - \sigma(\tilde{\tau}_-, X_{\tilde{\tau}})| + K\left[\int_{\tilde{\tau}_-}^{\tilde{\tau}} |b(t, X_t)|\, dt + \left|\int_{\tilde{\tau}_-}^{\tilde{\tau}} \gamma_t\, dW_t\right|\right]$$

$$\leq |\sigma(\tilde{\tau}, X_{\tilde{\tau}}) - \sigma(\tilde{\tau}_-, X_{\tilde{\tau}})| + C_0[(\tilde{\tau} - \tilde{\tau}_-) + \xi_T \gamma_T^* (\tilde{\tau} - \tilde{\tau}_-)^{1/3}].$$

Thus, one of the following three inequalities must hold true:

$$|\sigma(\tilde{\tau}, X_{\tilde{\tau}}) - \sigma(\tilde{\tau}_-, X_{\tilde{\tau}})| \geq \frac{T^\beta \gamma_T^*}{4},$$

$$(\tilde{\tau} - \tilde{\tau}_-)^{-1} \leq \frac{8C_0}{T^\beta \gamma_T^*},$$

$$(\tilde{\tau} - \tilde{\tau}_-)^{-1/3} \leq \frac{8C_0 \xi_T}{T^\beta}.$$

From the first inequality above, we have $\tilde{\tau} - \tilde{\tau}_- \geq \delta(\frac{T^\beta \gamma_T^*}{4})$. So in all three cases, we have

$$(5.7) \qquad (\tilde{\tau} - \tilde{\tau}_-)^{-1} \leq C\left[\delta\left(\frac{T^\beta \gamma_T^*}{4}\right)^{-1} + T^{-\beta}|\gamma_T^*|^{-1} + T^{-3\beta}\xi_T^3\right].$$

Similarly, if $\tau_+ < T$, then

$$(5.8) \qquad (\tilde{\tau}_+ - \tilde{\tau})^{-1} \leq C\left[\delta\left(\frac{T^\beta \gamma_T^*}{4}\right)^{-1} + T^{-\beta}|\gamma_T^*|^{-1} + T^{-3\beta}\xi_T^3\right].$$

Finally, if $\tau_- = 0$ and $\tau_+ = T$, then $(\tau_+ - \tau_-)^{-1} = T^{-1}$. Moreover, in this case we have $|\gamma_t| \geq (1 - \frac{T^\beta}{2})\gamma_T^*$ for all $t \in [0, T]$. In particular, we have

$$\left(1 - \frac{T^\beta}{2}\right)|\gamma_0| \leq \left(1 - \frac{T^\beta}{2}\right)\gamma_T^* \leq |\gamma_t| \leq \gamma_T^* \leq \left(1 - \frac{T^\beta}{2}\right)^{-1}|\gamma_0|.$$

Note that $\gamma_t$ is continuous, so all $\gamma_t$ have the same sign. Then one can easily prove that $|\gamma_t - \gamma_0| \leq T^\beta |\gamma_0|$, $\forall t$. That is, $\mathbb{1}_{G_T} = 1$.

Combining the above three cases, we know that the following inequality holds true:

$$(5.9) \quad (\tilde{\tau}_+ - \tilde{\tau}_-)^{-1} \leq C\left[\delta\left(\frac{T^\beta \gamma_T^*}{4}\right)^{-1} + T^{-\beta}|\gamma_T^*|^{-1} + T^{-3\beta}\xi_T^3 + T^{-1}\mathbb{1}_{G_T}\right],$$

which, combined with (5.6), obviously implies (5.5), and whence the theorem. $\square$

We note that if $\sigma$ is uniformly Hölder-$\alpha$ continuous in $t$, then $\delta(\varepsilon) \geq \varepsilon^{1/\alpha}$. The following result is a direct consequence of Lemma 5.6.

COROLLARY 5.7. *Assume* (A1), (A2$'$) *and* (A5). *Then*

$$|\Lambda_r^{t,x}|^{-1} \leq C|\gamma_r^{t,x,*}|^{-2}\bigg[(r-t)^{-\beta/\alpha}|\gamma_r^{t,x,*}|^{-1/\alpha}$$
$$+ (r-t)^{-\beta}|\gamma_r^{t,x,*}|^{-1} + \frac{|\xi_r^{t,x}|^3}{(r-t)^{3\beta}} + \frac{\mathbb{1}_{G_r^{t,x}}}{r-t}\bigg].$$

For notational convenience, in the sequel we denote $E\{\xi; A\} \stackrel{\triangle}{=} E\{\xi \mathbb{1}_A\}$ for a random variable $\xi$ and an event $A$.

LEMMA 5.8. *Assume* (A1), (A2$'$), (A5) *and* $T \leq T_0$. *Recall* (4.12) *and denote* $p \stackrel{\triangle}{=} \frac{24\alpha}{1+22\alpha} > 1$, *where* $\alpha$ *is as in* (A2$'$). *Then for any* $c > 0$, *there exists a constant* $C$, *depending only on* $K, T_0, \alpha$ *and* $c$, *such that:*

(i) $E^{0,x}\{|N_t|^p; |\Delta X_t| \geq ct\} \leq Ct^{-p}$;
(ii) $\int_0^T E^{0,x}\{|N_t|^p; |\Delta X_t| \geq ct\}^{1/p} dt \leq C$.

PROOF. Again we omit the superscript $^{0,x}$. Denote

$$\beta \stackrel{\triangle}{=} \frac{2\alpha - 1}{12\alpha}; \qquad \mu \stackrel{\triangle}{=} \frac{1 + 10\alpha}{24\alpha}$$

and

$$I_t \stackrel{\triangle}{=} E\{|N_t|^p; |\Delta X_t| \geq ct\}^{1/p}; \qquad I \stackrel{\triangle}{=} \int_0^T I_t\, dt.$$

Recalling (4.7), we have

$$I_t^p \leq CE\bigg\{\frac{1}{\Lambda_t^p}\bigg[\bigg|\int_0^t \gamma_s \nabla X_s\, dW_s\bigg|^p + \bigg|\int_0^t \frac{\Lambda_s}{\Lambda_t}\sigma_x(s, X_s)\gamma_s \nabla X_s\, ds\bigg|^p\bigg];$$

(5.10)
$$|\Delta X_t| \geq ct\bigg\}$$

$$\leq CE\bigg\{\frac{1}{\Lambda_t^p}|\gamma_t^* \nabla X_t^*|^p[\tilde{\xi}_t^p + t^p]; |\Delta X_t| \geq ct\bigg\},$$

where

$$\nabla X_t^* \stackrel{\triangle}{=} \sup_{0 \leq s \leq t} |\nabla X_s|; \qquad \tilde{\xi}_t \stackrel{\triangle}{=} \bigg[\sup_{0 \leq s \leq t} |\gamma_s \nabla X_s|\bigg]^{-1} \int_0^t \gamma_s \nabla X_s\, dW_s.$$

Now applying Corollary 5.7, one has

(5.11) $$I_t^p \leq C[I_{1,t}^p + I_{2,t}^p],$$

where

$$I_{1,t}^p \triangleq E\{|\gamma_t^*|^{-p}|\nabla X_t^*|^p(t^p + \tilde{\xi}_t^p)[|t^\beta \gamma_t^*|^{-p/\alpha} + t^{-p\beta}|\gamma_t^*|^{-p} + t^{-3p\beta}\xi_t^{3p}];$$
$$|\Delta X_t| \geq ct\};$$

$$I_{2,t}^p \triangleq t^{-p}E\{|\gamma_t^*|^{-p}|\nabla X_t^*|^p[t^p + \tilde{\xi}_t^p]; G_t, |\Delta X_t| \geq ct\}.$$

We estimate $I_{1,t}^p$ first. By Lemma 2.3, it holds that $E\{|\nabla X_t^*|^q\} \leq C_q, \forall q \geq 1$. Recalling (5.1), we have

$$(5.12) \quad I_{1,t}^p \leq Ct^{p/2}E\{t^{-2p\beta/\alpha}|\gamma_t^*|^{-2p(1+1/\alpha)} + t^{-2p\beta}|\gamma_t^*|^{-4p}$$
$$+ t^{-6p\beta}|\gamma_t^*|^{-2p}\xi_t^{6p}; |\Delta X_t| \geq ct\}^{1/2}.$$

Now we recall (4.14). Note that $\Delta X_t = L_t^{-1}M_t$ and both $L_t$ and $L_t^{-1}$ are bounded, then by (5.1), one can easily show that, for all $q \geq 2$,

$$E\{|\gamma_t^*|^{-q}; |\Delta X_t| \geq ct\} \leq Ct^{-q}E\{|\gamma_t^*|^{-q}|\Delta X_t|^q\}$$
$$\leq Ct^{-q}E\{|(L\gamma)_t^*|^{-q}|M_t|^q\} \leq Ct^{-q/2}.$$

Thus, by (5.12) and applying Lemma 5.5, we have

$$I_{1,t}^p \leq Ct^{p/2}[t^{-(p/(2\alpha))(2\beta+\alpha+1)} + t^{-p(\beta+1)} + t^{-3\beta}] \leq Ct^{-(p/2\alpha)(2\beta+1)},$$

thanks to the fact that $\beta < \frac{1}{4}$. Since $\alpha > \frac{1}{2}$, we have $2\beta + 1 < 2\alpha$ and therefore,

$$(5.13) \quad I_{1,t}^p \leq Ct^{-p}; \qquad \int_0^T I_{1,t}\,dt \leq C\int_0^T t^{-(2\beta+1)/2\alpha}\,dt \leq C.$$

It remains to estimate $I_{2,t}$. To this end, we recall that

$$G_t = \left\{\sup_{0 \leq s \leq t} |\gamma_s - \gamma_0| \leq |\gamma_0|t^\beta\right\}.$$

If $\gamma_0 = 0$, then in $G_t$, we have $\gamma_s = 0$ for all $s \in [0,t]$ and, thus, $\Delta X_t = 0$. That is, if $\gamma_0 = 0$, then $P\{G_t, \Delta X_t \geq ct\} = 0$, which implies that $I_{2,t} = 0$. Now we assume $\gamma_0 \neq 0$. Without loss of generality, we assume $\gamma_0 > 0$. We note again that in the sequel the constants $C$ and $c$ may vary from line to line, but they are independent of $\frac{1}{\gamma_0}$. One can check directly that

$$\{G_t, |\Delta X_t| \geq ct\} \subset \{G_t, |M_t| \geq ct\}$$
$$\subset \left\{\gamma_0|W_t| \geq \frac{ct}{2}\right\} \cup \left\{G_t, \left|\int_0^t [L_s\gamma_s - \gamma_0]\,dW_s\right| \geq \frac{ct}{2}\right\}.$$

Denote

$$\lambda_s^t \triangleq \frac{L_s\gamma_s - \gamma_0}{\gamma_0 t^\beta}; \qquad \xi \triangleq \left[\sup_{0 \leq s \leq t}|\gamma_s^t|\right]^{-1}\left|\int_0^t \lambda_s^t\,dW_s\right|.$$

Note that $\beta < 1$. In $G_t$ and for $\forall\, s \leq t$, we have

$$|L_s\gamma_s - \gamma_0| \leq |L_s - 1|\gamma_0 + L_s|\gamma_s - \gamma_0| \leq C[\gamma_0 s + \gamma_0 t^\beta] \leq C\gamma_0 t^\beta,$$

which implies that $|\lambda_s^t| \leq C$ in $G_t$. Thus,

$$\left\{G_t, \left|\int_0^t [L_s\gamma_s - \gamma_0]\, dW_s\right| \geq ct\right\} \subset \left\{G_t, \gamma_0 t^\beta \left|\int_0^t \lambda_s^t\, dW_s\right| \geq ct\right\}$$
$$\subset \{\gamma_0 t^\beta \xi \geq ct\}.$$

Therefore,

$$\{G_t, |\Delta X_t| \geq ct\} \subset \{\gamma_0 |W_t| \geq ct\} \cup \{\gamma_0 t^\beta \xi \geq ct\}.$$

Now by the definition of $I_{2,t}$, we have

(5.14) $$I_{2,t}^p \leq I_{3,t}^p + I_{4,t}^p,$$

where

$$I_{3,t}^p \triangleq t^{-p} E\{|\gamma_t^*|^{-p}|\nabla X_t^*|^p [t^p + \tilde{\xi}_t^p]; \gamma_0|W_t| \geq ct\};$$
$$I_{4,t}^p \triangleq t^{-p} E\{|\gamma_t^*|^{-p}|\nabla X_t^*|^p [t^p + \tilde{\xi}_t^p]; \gamma_0 t^\beta \xi \geq ct\}.$$

Note that

(5.15)
$$I_{3,t}^p \leq t^{-p} E\left\{|\gamma_t^*|^{-p}|\nabla X_t^*|^p [t^p + \tilde{\xi}_t^p]\left(\frac{\gamma_0 |W_t|}{ct}\right)^p\right\}$$
$$\leq Ct^{-2p} E\{|\nabla X_t^*|^p [t^p + \tilde{\xi}_t^p]|W_t|^p\}$$
$$\leq Ct^{-2p} E\{|\nabla X_t^*|^{3p}\}^{1/3} E\{t^{3p} + \tilde{\xi}_t^{3p}\}^{1/3} E\{|W_t|^{3p}\}^{1/3}$$
$$\leq Ct^{-2p} t^{p/2} t^{p/2} = Ct^{-p},$$

where the last inequality is thanks to Lemma 2.3 and (5.1). Similarly, applying Lemma 5.5, one can prove

(5.16) $$I_{4,t}^p \leq Ct^{p(\beta-1)}.$$

Since $\alpha > \frac{1}{2}$, one has $\beta > 0$. Then by (5.14), (5.15) and (5.16), we have

$$I_{2,t}^p \leq Ct^{-p},$$

which, combined with (5.13), proves (i).

To prove (ii), first by (5.16), we have

(5.17) $$\int_0^T I_{4,t}\, dt \leq C \int_0^T t^{\beta-1}\, dt = C < \infty.$$

Moreover, by the definition of $I_{3,t}$ and applying Lemma 2.3 and (5.1), again we have

$$\int_0^T I_{3,t}\,dt \le \int_0^T (\gamma_0 t)^{-1} E\{|\nabla X_t^*|^p [t^p + \tilde{\xi}_t^p]; \gamma_0 |W_t| \ge ct\}^{1/p}\,dt$$

$$= \left[\int_0^{\gamma_0^2} + \int_{\gamma_0^2}^T\right] (\gamma_0 t)^{-1} E\{|\nabla X_t^*|^p [t^p + \tilde{\xi}_t^p]; \gamma_0 |W_t| \ge ct\}^{1/p}\,dt$$

$$\le \int_0^{\gamma_0^2} \frac{1}{\gamma_0 t} E\{|\nabla X_t^*|^p [t^p + \tilde{\xi}_t^p]\}^{1/p}\,dt$$

$$+ \int_{\gamma_0^2}^T \frac{1}{\gamma_0 t} E\left\{|\nabla X_t^*|^p [t^p + \tilde{\xi}_t^p]\left(\frac{\gamma_0 |W_t|}{ct}\right)^{2p}\right\}^{1/p}\,dt$$

$$\le C \int_0^{\gamma_0^2} (\gamma_0 t)^{-1} t^{1/2}\,dt + C \int_{\gamma_0^2}^T (\gamma_0 t)^{-1} t^{1/2} \gamma_0^2 t^{-1}\,dt$$

$$= C\gamma_0^{-1} t^{1/2}\big|_{t=0}^{t=\gamma_0^2} - C\gamma_0 t^{-1/2}\big|_{t=\gamma_0^2}^{t=T} \le C < \infty,$$

which, combined with (5.17) and (5.13), proves (ii). □

5.3. *A localizing result.* In this section we prove a localizing version of Theorem 3.1, which will play a very important role in the proof of Theorem 5.1. To this end, we first introduce a notion called "$\varepsilon$-neighbor."

DEFINITION 5.9. Fix $K$ as an upper bound of $|b_x|$. For any $(t_0, x_0) \in [0, T) \times \mathbb{R}$ and any $\varepsilon > 0$, the $\varepsilon$-neighbor of $(t_0, x_0)$ is the set

$$D_\varepsilon(t_0, x_0) \triangleq \{(t, x) : t_0 \le t \le T, |x - \eta_t^{t_0, x_0}| \le \varepsilon[e^{K(t-t_0)} - 1]\}.$$

The following lemma gives a basic property of $\varepsilon$-neighbors.

LEMMA 5.10. *If $(t_1, x_1) \in D_\varepsilon(t_0, x_0)$, then $D_\varepsilon(t_1, x_1) \subset D_\varepsilon(t_0, x_0)$.*

PROOF. Assume $(t_1, x_1) \in D_\varepsilon(t_0, x_0)$ and $(t_2, x_2) \in D_\varepsilon(t_1, x_1)$. Denote $\eta^i \triangleq \eta^{t_i, x_i}$ for $i = 0, 1$. Then

$$|x_1 - \eta_{t_1}^0| \le \varepsilon[e^{K(t_1 - t_0)} - 1]; \qquad |x_2 - \eta_{t_2}^1| \le \varepsilon[e^{K(t_2 - t_1)} - 1].$$

Note that, for $t \ge t_1$,

$$\eta_t^0 = \eta_{t_1}^0 + \int_{t_1}^t b(s, \eta_s^0)\,ds; \qquad \eta_t^1 = x_1 + \int_{t_1}^t b(s, \eta_s^1)\,ds.$$

Denote

$$\Delta \eta_t \triangleq \eta_t^1 - \eta_t^0; \qquad \alpha_t \triangleq \frac{b(t, \eta_t^1) - b(t, \eta_t^0)}{\Delta \eta_t}.$$

Then $|\alpha_t| \leq K$, and

$$\Delta \eta_t = x_1 - \eta_{t_1}^0 + \int_{t_1}^t \alpha_s \Delta \eta_s \, ds.$$

Thus,

$$\Delta \eta_t = (x_1 - \eta_{t_1}^0) \exp\left(\int_{t_1}^t \alpha_s \, ds\right).$$

Therefore,

$$|\Delta \eta_{t_2}| \leq |x_1 - \eta_{t_1}^0| \exp\left(\int_{t_1}^{t_2} \alpha_s \, ds\right) \leq \varepsilon [e^{K(t_1 - t_0)} - 1] e^{K(t_2 - t_1)}.$$

Then

$$|x_2 - \eta_{t_2}^0| \leq |x_2 - \eta_{t_2}^1| + |\Delta \eta_{t_2}|$$
$$\leq \varepsilon [e^{K(t_2 - t_1)} - 1] + \varepsilon [e^{K(t_2 - t_0)} - e^{K(t_2 - t_1)}] = \varepsilon [e^{K(t_2 - t_0)} - 1],$$

which proves that $(t_2, x_2) \in D_\varepsilon(t_0, x_0)$. □

The following lemma is the key part for the proof of Theorem 5.1.

LEMMA 5.11. *Assume that all the conditions in Theorem 5.1 hold true. Assume further that $T \leq T_0$ and that, for some $(t_0, x_0) \in [0, T) \times \mathbb{R}$ and some constants $\varepsilon_1, \varepsilon_2 > 0$, $|u_x(t,x)| \leq K_0$ for $\forall (t,x) \in D_{\varepsilon_2}(t_0, x_0) \cap ([T - \varepsilon_1, T] \times \mathbb{R})$. Then*

$$|u_x(t_0, x_0)| \leq C[K_0 + \psi(x_0)],$$

*where $C$ depends on $K, \alpha, \varepsilon_1, \varepsilon_2$ and $T_0$, but* not *on the upper bound of $g'$.*

PROOF. Without loss of generality, we assume $t_0 = 0$ and omit the superscript $^{0,x_0}$. Denote $D_{\varepsilon_2} \triangleq D_{\varepsilon_2}(0, x_0)$. For $t \in [0, T]$, denote

$$A_t \triangleq \sup_{\{x \,:\, (t,x) \in D_{\varepsilon_2}\}} |u_x(t, x)|.$$

Applying Theorem 3.1, we have $\sup_{0 \leq t \leq T} A_t < \infty$. Moreover, by assumption, we have

(5.18) $$A_t \leq K_0 \qquad \forall t \in [T - \varepsilon_1, T].$$

We claim that

(5.19) $$A_t \leq C\left[K_0 + \sup_{\{x \,:\, (t,x) \in D_{\varepsilon_2}\}} \psi(x) + \int_t^T A_s \, ds\right] \qquad \forall t \in [0, T - \varepsilon_1].$$

For $(t,x) \in D_{\varepsilon_2}$, we have
$$|x| \leq |x - \eta_t| + |\eta_t - x_0| + |x_0| \leq \varepsilon_2[e^{Kt} - 1] + Ct + |x_0| \leq |x_0| + C.$$
Thus,
$$\psi(x) \leq C[\psi(x_0) + \psi(C)] \leq C\psi(x_0),$$
which, combined with (5.19) and (5.18), implies that
$$A_t \leq C\left[K_0 + \psi(x_0) + \int_t^T A_s \, ds\right]$$
$$\leq C\left[K_0 + \psi(x_0) + \int_t^{T-\varepsilon_1} A_s \, ds\right] \qquad \forall t \in [0, T - \varepsilon_1].$$
Then the lemma follows the Gronwall inequality.

It remains to prove (5.19). By Lemma 5.10, it suffices to prove it at $t = 0$. In this case, (5.19) becomes
$$(5.20) \qquad |u_x(0, x_0)| \leq C\left[K_0 + \psi(x_0) + \int_0^T A_t \, dt\right].$$
Note that
$$u(0, x) = E^x\left\{g(X_T) + \int_0^T f(r, X_r, Y_r) \, dr\right\}.$$
Thus,
$$u_x(0, x_0) = E\left\{g'(X_T)\nabla X_T + \int_0^T [f_x \nabla X_t + f_y \nabla Y_t] \, dt\right\}.$$
For any $t \in (0, T]$, let $\varphi(t, x)$ be a smooth function of $x$ satisfying that
$$\varphi(t, x) = \begin{cases} 1, & \text{if } |x - \eta_t| \geq \varepsilon_2[e^{Kt} - 1], \\ 0, & \text{if } |x - \eta_t| \leq \frac{\varepsilon_2}{2}[e^{Kt} - 1], \end{cases}$$
$$|\varphi(t, x)| \leq 1; \qquad |\varphi_x(t, x)| \leq \frac{C}{t}.$$
Then we have
$$(5.21) \qquad \begin{aligned} u_x(0, x_0) = E\Big\{ & g'(X_T)\nabla X_T \varphi(T, X_T) \\ & + \int_0^T [f_x \nabla X_t + f_y \nabla Y_t]\varphi(t, X_t) \, dt \\ & + g'(X_T)\nabla X_T(1 - \varphi(T, X_T)) \\ & + \int_0^T [f_x \nabla X_t + f_y \nabla Y_t](1 - \varphi(t, X_t)) \, dt \Big\}. \end{aligned}$$

Since $\nabla$ can be considered as the differential operator with respect to $x$, by the chain rule, one can easily get that

$$g'(X_T)\nabla X_T \varphi(T, X_T) = [\nabla g(X_T)]\varphi(T, X_T)$$
$$= \nabla(g(X_T)\varphi(T, X_T)) - g(X_T)\varphi_x(T, X_T)\nabla X_T.$$

Note that when $\varphi(T, X_T) \neq 0$, one has $|\Delta X_T| \geq \frac{\varepsilon_2}{2}(e^{KT} - 1)$, which implies that $(0, x_0) \in \Gamma^0$ and, thus, $\Lambda_T > 0$. Then following the arguments in Section 4.3, one can easily prove that

$$E\{g'(X_T)\nabla X_T \varphi(T, X_T)\} = E\{g(X_T)[\varphi(T, X_T)N_T^0 - \varphi_x(T, X_T)\nabla X_T]\}.$$

Similarly, we have

$$E\{[f_x \nabla X_t + f_y \nabla Y_t]\varphi(t, X_t)\} = E\{f(t, X_t, Y_t)[\varphi(t, X_t)N_t^0 - \varphi_x(t, X_t)\nabla X_t]\}.$$

So one can rewrite (5.21) as

$$u_x(0, x_0) = I_1 + I_2 + I_3,$$

where

$$I_1 \triangleq E\{g(X_T)[\varphi(T, X_T)N_T - \varphi_x(T, X_T)\nabla X_T]\}$$

$$I_2 \triangleq E\left\{\int_0^T f(t, X_t, Y_t)[\varphi(t, X_t)N_t - \varphi_x(t, X_t)\nabla X_t]\, dt\right\}$$

$$I_3 \triangleq E\left\{g'(X_T)\nabla X_T(1 - \varphi(T, X_T))\right.$$
$$\left. + \int_0^T [f_x \nabla X_t + f_y u_x \nabla X_t](1 - \varphi(t, X_t))\, dt\right\}.$$

We shall estimate $I_1 - I_3$ separately. First, it is obvious that

$$(5.22) \qquad |I_3| \leq C\left[A_T + \int_0^T (1 + A_t)\, dt\right] \leq C\left[K_0 + \int_0^T A_t\, dt\right].$$

Denote $p \triangleq \frac{24\alpha}{1+22\alpha} > 1$ as in Lemma 5.8. Let $q$ be the conjugate of $p$. Applying Lemma 5.8, we have

$$E\{|g(X_T)\varphi(T, X_T)N_T|\} \leq E\left\{|g(X_T)N_T|; |\Delta X_T| \geq \frac{\varepsilon_2}{2}[e^{KT} - 1]\right\}$$

$$\leq E\{|g(X_T)N_T|; |\Delta X_T| \geq cT\}$$

(5.23)

$$\leq \|g(X_T)\|_q [E\{|N_T|^p; |\Delta X_T| \geq cT\}]^{1/p}$$

$$\leq C\psi(x_0)T^{-1} \leq C\psi(x_0),$$

where the last inequality is thanks to the fact that $T \geq \varepsilon_1$. Moreover,

$$E\{|g(X_T)\varphi_x(T, X_T)\nabla X_T|\} \leq \frac{C}{T}E\{|g(X_T)\nabla X_T|\} \leq C\psi(x_0),$$

which, combined with (5.23), implies that

(5.24) $$|I_1| \leq C\psi(x_0).$$

It remains to estimate $I_2$. By the arguments in Section 4.3, we have

$$E\{\varphi_x(t, X_t)\nabla X_t\} = E\{\varphi(t, X_t)N_t\} \qquad \forall t \in (0, T].$$

Thus,

(5.25) $$I_2 = E\bigg\{\int_0^T [f(t, X_t, u(t, X_t)) - f(t, \eta_t, u(t, \eta_t))]$$
$$\times [\varphi(t, X_t)N_t - \varphi_x(t, X_t)\nabla X_t]\, dt\bigg\}.$$

By (2.10),

(5.26) $$|u(t,x)| = |Y_t^{t,x}| \leq C\psi(x),$$

which implies that

$$|f(t, \eta_t, u(t, \eta_t))| \leq C[1 + |\eta_t| + |u(t, \eta_t)|] \leq C\psi(\eta_t);$$
$$|f(t, X_t, u(t, X_t))| \leq C\psi(X_t).$$

Then by (2.1) and applying Lemma 5.8, again we have

(5.27) $$\int_0^T E\{|[f(t, X_t, u(t, X_t)) - f(t, \eta_t, u(t, \eta_t))]\varphi(t, X_t)N_t|\}\, dt$$
$$\leq C\psi(x_0)\int_0^T E\{|N_t|^p; |\Delta X_t| > ct\}^{1/p}\, dt \leq C\psi(x_0).$$

Moreover, noting that $\varphi_x(t, X_t) = 0$ when $|\Delta X_t| \geq \varepsilon_2[e^{Kt} - 1]$, we have

$$\int_0^T E\{|[f(t, X_t, u(t, X_t)) - f(t, \eta_t, u(t, \eta_t))]\varphi_x(t, X_t)\nabla X_t|\}\, dt$$
$$\leq C\int_0^T E\bigg\{\bigg|[f(t, X_t, u(t, X_t)) - f(t, \eta_t, u(t, \eta_t))]\frac{1}{t}\nabla X_t\bigg|;$$
$$|\Delta X_t| \leq \varepsilon_2[e^{Kt} - 1]\bigg\}\, dt$$
$$\leq C\int_0^T E\{|\Delta X_t|(1 + A_t)\nabla X_t|; |\Delta X_t| \leq \varepsilon_2[e^{Kt} - 1]\}\frac{dt}{t}$$
$$\leq C\bigg[1 + \int_0^T A_t\, dt\bigg],$$

which, combined with (5.27) and (5.25), implies that

(5.28) $$|I_2| \leq C\bigg[\psi(x_0) + \int_0^T A_t\, dt\bigg].$$

Now combining (5.22), (5.24) and (5.28), we prove (5.20), and hence the lemma. □

5.4. *Proof of Theorem* 5.1. We first prove a simple case.

LEMMA 5.12. *Assume that all the conditions in Theorem* 5.1 *hold true, and that* $|\sigma(t_0, x_0)| \geq \frac{1}{n}$. *Then the result of Theorem* 5.1 *holds true at* $(t_0, x_0)$.

PROOF. Without loss of generality, we assume $t_0 = 0$ and omit the superscript $^{0,x_0}$ as before. Then $\gamma_t^* \geq |\sigma(0, x_0)| \geq \frac{1}{n}$ for $\forall t \geq 0$. Choosing $\beta = 0$ in Lemma 5.6, we get

$$|\Lambda_t^{-1}| \leq C_n[1 + |\xi_t|^3 + t^{-1}]. \tag{5.29}$$

Then one can easily get

$$E\{|N_t^0|^2|\} \leq \frac{C_n}{t} \qquad \forall t > 0. \tag{5.30}$$

By Theorem 3.1, we know $u_x(0, x_0)$ exists. By the same arguments as in Section 4.3 or as in the proof of Lemma 5.11, one can show that

$$u_x(0, x_0) = E\Big\{g(X_T)N_T^0 + \int_0^T f(t, X_t, Y_t)N_t^0 \, dt\Big\},$$

which, together with (2.1) and (5.26), implies that

$$|u_x(0, x_0)| \leq C_n\Big[\frac{\psi(x_0)}{\sqrt{T}} + \int_0^T \frac{\psi(x_0)}{\sqrt{t}} \, dt\Big] \leq \frac{C_n \psi(x_0)}{\sqrt{T}}.$$

That proves the lemma. □

PROOF OF THEOREM 5.1. Fix $(t_0, x_0) \in \Gamma^n$. Without loss of generality, we assume $t_0 = 0$ again. First if $|\sigma(0, x_0)| \geq \frac{1}{2n}$, then by Lemma 5.12, the results hold true.

Now we assume $|\sigma(0, x_0)| < \frac{1}{2n}$. Since $(0, x_0) \in \Gamma^n$, there exists $t_1 \in (0, T]$ such that $|\sigma(t_1, \eta_{t_1})| \geq \frac{1}{n}$. Since $\sigma, b$ are Lipschitz continuous in $x$ and $\sigma$ is Hölder-$\alpha$ continuous in $t$, there exist constants $\varepsilon_1, \varepsilon_2 > 0$, depending only on $K, n$ and $\alpha$ such that for any $(t, x) \in \bar{D}_{\varepsilon_1, \varepsilon_2}(t_1, \eta_{t_1}) \triangleq \{(t, x) : t \in [t_1 - \varepsilon_1, t_1], |x - \eta_t| \leq \varepsilon_2\}$, it holds that $|\sigma(t, x) - \sigma(t_1, \eta_{t_1})| \leq \frac{1}{2n}$. Thus, $|\sigma(t, x)| \geq \frac{1}{2n}$ for any $(t, x) \in \bar{D}_{\varepsilon_1, \varepsilon_2}(t_1, \eta_{t_1})$. We note that this also implies that $t_1 > \varepsilon_1$. Now we choose $t_2 \triangleq t_1 - \frac{\varepsilon_1}{2}$, $\delta_2 \triangleq \frac{\varepsilon_1}{2}$. Obviously, there exists $\delta_3 > 0$, depending only on $K, n$ and $\alpha$ such that $D_{\delta_3}(0, x_0) \cap ([t_2 - \delta_2, t_2] \times \mathbb{R}) \subset \bar{D}_{\varepsilon_1, \varepsilon_2}(t_1, \eta_{t_1})$. Now for $\forall (t, x) \in D_{\delta_3}(0, x_0) \cap ([t_2 - \delta_2, t_2] \times \mathbb{R})$, by Lemma 5.12, we have

$$|u_x(t, x)| \leq \frac{C_n \psi(x)}{\sqrt{T-t}} \leq \frac{C_n \psi(x)}{\sqrt{t_1 - t_2}} = \frac{C_n \psi(x)}{\sqrt{\delta_2}} = C_n \psi(x) \leq C_n \psi(x_0),$$

where the last inequality is thanks to the fact that $(t,x) \in D_{\delta_3}(0,x_0)$. Then applying Lemma 5.11, we get $|u_x(0,x_0)| \leq C_n[C_n\psi(x_0) + \psi(x_0)] = C_n\psi(x_0)$. $\square$

5.5. *Proof of Theorem* 5.2. We first prove the theorem under (A3$'$) instead of (A3$''$).

LEMMA 5.13. *Assume* (A1), (A2$'$), (A3$'$), (A4) *and* (A5). *Then all the results in Theorem* 5.2 *hold true*.

PROOF. First if $f \in C_b^{0,1}, g \in C_b^1$, then by Theorem 3.1, we know $u \in C^{0,1}([0,T] \times \mathbb{R}) \subset C^{0,1}(\Gamma)$ and $Z_t = (u_x\sigma)(t,X_t)$. Since $\sigma(t,X_t) = 0$ for $t > \tau$, so (iii) holds true. (ii) is due to Theorem 5.1. Finally, for $(t,x) \in \Gamma$, we have $|\sigma(t,x)| \geq \frac{1}{n}$ for some $n$. Then the representation formula in (i) follows the proof of Lemma 5.12.

For the general case, one can easily prove the theorem by following the approximating arguments in the proof of Theorem 4.4 (or Theorem 3.2). $\square$

We now assume only (A3$''$). In this case (1.2) can be rewritten as

$$
\begin{aligned}
&u_t + \tfrac{1}{2}\sigma^2 u_{xx} + \tilde{b}u_x + f_1(t,x,u) = 0; \\
&u(T,x) = g(x),
\end{aligned}
\tag{5.31}
$$

where $\tilde{b} \triangleq b + f_2\sigma$. Recall that $(X,Y,Z)$ is the solution to FBSDE (1.1). Define

$$\tilde{\eta}_t = x + \int_0^t \tilde{b}(r,\tilde{\eta}_r)\,dr; \qquad \tilde{\Gamma}^0 \triangleq \Big\{(t,x) \colon \sup_{t\leq s\leq T}|\sigma(s,\tilde{\eta}_s^{t,x})| > 0\Big\};$$

$$\tilde{\tau} \triangleq \inf\{t \colon (t,X_t) \notin \tilde{\Gamma}^0\}; \qquad \tilde{\Gamma}^n \triangleq \Big\{(t,x) \colon \sup_{t\leq s\leq T}|\sigma(s,\tilde{\eta}_s^{t,x})| \geq \frac{1}{n}\Big\};$$

$$d\tilde{W}_t \triangleq dW_t - f_2(t,X_t)\,dt.$$

Then $\tilde{W}$ is a Brownian motion under another probability and one can rewrite (1.1) as

$$
\begin{aligned}
X_t &= x + \int_0^t \tilde{b}(r,X_r)\,dr + \int_0^t \sigma(r,X_r)\,d\tilde{W}_r; \\
Y_t &= g(X_T) + \int_t^T f_1(r,X_r,Y_r)\,dr - \int_t^T Z_r\,d\tilde{W}_r.
\end{aligned}
\tag{5.32}
$$

LEMMA 5.14. *Assume* (A1) *and* (A3$''$). *Then*:

(i) $\tilde{\Gamma}^0 = \Gamma^0$, $\tilde{\tau} = \tau$.

(ii) *There exist a constant $C > 1$, depending only on $K$, such that $\Gamma^n \subset \tilde{\Gamma}^{Cn}$ and $\tilde{\Gamma}^n \subset \Gamma^{Cn}$, for any $n \geq 1$.*

PROOF. (i) If $(t, x) \notin \Gamma^0$, then $\sigma(s, \eta_s^{t,x}) = 0$ for $\forall t \leq s \leq T$. So it holds that $\tilde{b}(s, \eta_s^{t,x}) = b(s, \eta_s^{t,x})$. Since $\eta$ is the solution to (4.2), we have $\eta_s^{t,x} = x + \int_t^s \tilde{b}(r, \eta_r^{t,x}) \, dr$, which implies that $\tilde{\eta}_s^{t,x} = \eta_s^{t,x}$. Therefore, $\sigma(s, \tilde{\eta}_s^{t,x}) = \sigma(s, \eta_s^{t,x}) = 0$, and, hence, $(t, x) \notin \tilde{\Gamma}^0$. That is, $\tilde{\Gamma}^0 \subset \Gamma^0$. Similarly, one can prove that $\Gamma^0 \subset \tilde{\Gamma}^0$. So $\tilde{\Gamma}^0 = \Gamma^0$. Then $\tilde{\tau} = \tau$ follows immediately.

(ii) Fix $(t, x)$. Denote $\Delta \eta_s \triangleq \tilde{\eta}_s^{t,x} - \eta_s^{t,x}$. Then

$$\Delta \eta_s = \int_t^s [\beta_r \Delta \eta_r + (f_2 \sigma)(r, \tilde{\eta}_r^{t,x})] \, dr,$$

where $\beta_r \triangleq \frac{b(r, \tilde{\eta}_r^{t,x}) - b(r, \eta_r^{t,x})}{\Delta \eta_r}$ is bounded. Thus, we have

$$\sup_{t \leq s \leq T} |\Delta \eta_s| \leq C \sup_{t \leq s \leq T} |\sigma(s, \tilde{\eta}_s^{t,x})|.$$

Note that

$$|\sigma(s, \eta_s^{t,x}) - \sigma(s, \tilde{\eta}_s^{t,x})| \leq C |\Delta \eta_s|.$$

Then $\sup_{t \leq s \leq T} |\sigma(s, \eta_s^{t,x})| \leq C \sup_{t \leq s \leq T} |\sigma(s, \tilde{\eta}_s^{t,x})|$, and therefore, $\Gamma^n \subset \tilde{\Gamma}^{Cn}$.

On the other hand, since $b = \tilde{b} - f_2 \sigma$, one can similarly show that $\tilde{\Gamma}^n \subset \Gamma^{Cn}$.

□

PROOF OF THEOREM 5.2. (iii) is a direct consequence of Lemmas 5.13(iii) and 5.14(i). As to (ii), for any $(t, x) \in \Gamma^n$, by Lemma 5.14(ii), $(t, x) \in \tilde{\Gamma}^{Cn}$ for some constant $C$. Applying Lemma 5.13, we get $u_x(t, x) \leq \frac{C_n \psi(x)}{\sqrt{T-t}}$.

It remains to prove (i). By Lemma 5.13 and (5.31), we have $u \in C^{0,1}(\Gamma)$. To prove the representation theorem, by standard approximating arguments, it suffices to show that the integral in the right-hand side of the formula converges. To this end, we fix an $(t_0, x_0) \in \Gamma$. That is, $|\sigma(t_0, x_0)| > 0$. Assume $|\sigma(t_0, x_0)| \geq \frac{1}{n}$ for some $n$. Without loss of generality, we assume again that $t_0 = 0$ and that $\sigma(0, x_0) \geq \frac{1}{n}$. Recall the proof of Lemma 5.12. By (5.30) and noting that $|Y_t| = |u(t, X_t)| \leq C \psi(X_t)$, one can easily show that

$$(5.33) \quad E\left\{|g(X_T) N_T^0| + \int_0^T |f_1(t, X_t, Y_t) N_t^0| \, dt \right\} \leq C_{n,T} \psi(x_0) < \infty.$$

Here $C_{n,T}$ may depend on $T^{-1}$ as well.

We finally show that

$$(5.34) \quad E\left\{ \int_0^T |f_2(t, X_t) Z_t N_t^0| \, dt \right\} \leq C_{n,T} \psi(x_0) < \infty.$$

To this end, we define
$$\tau_1 \triangleq \inf\left\{t : \sigma(t, X_t) = \frac{1}{2n}\right\} \wedge T.$$

Then for $t \leq \tau_1$, it holds that $\sigma(t, X_t) \geq \frac{1}{2n}$. Thus, by (ii) and (iii), we have $|Z_t| = |(u_x \sigma)(t, X_t)| \leq \frac{C_n \psi(X_t)}{\sqrt{T-t}}$. Therefore,

$$\begin{aligned}
E&\left\{\int_0^T |f_2(t, X_t) Z_t N_t^0| \, dt\right\} \\
&\leq CE\left\{\left[\int_0^{\tau_1} + \int_{\tau_1}^T\right] |Z_t N_t^0| \, dt\right\} \\
&\leq C_n E\left\{\int_0^{\tau_1} \frac{1}{\sqrt{T-t}} \psi(X_t) |N_t^0| \, dt\right\} + CE\left\{\int_{\tau_1}^T \frac{1}{\sqrt{t}} |Z_t \sqrt{t} N_t^0| \, dt\right\} \\
&\leq C_n \int_0^T \frac{\psi(x_0)}{\sqrt{t(T-t)}} \, dt \\
&\quad + CE\left\{\int_0^T |Z_t|^2 \, dt\right\}^{1/2} E\left\{\int_0^T t^2 |N_t^0|^4 \, dt\right\}^{1/4} E\left\{\int_{\tau_1}^T \frac{1}{t^2} \, dt\right\}^{1/4} \\
&\leq C_n \psi(x_0) + C_n \psi(x_0) E\{\tau_1^{-1}\}^{1/4}.
\end{aligned}$$
(5.35)

We follow the arguments in Lemma 5.6 to estimate $\tau_1^{-1}$. If $\tau_1 = T$, then $\tau_1^{-1} = \frac{1}{T}$. Now we assume $\tau_1 < T$, then $\sigma(\tau_1, X_{\tau_1}) = \frac{1}{2n}$. Note that $\sigma(0, x_0) \geq \frac{1}{n}$. Thus,

$$\frac{1}{2n} \leq \sigma(0, x_0) - \sigma(\tau_1, X_{\tau_1}) \leq |\sigma(0, x_0) - \sigma(\tau_1, x_0)| + |\sigma(\tau_1, x_0) - \sigma(\tau_1, X_{\tau_1})|$$

$$\leq C[\tau_1^\alpha + |X_{\tau_1} - x_0|] \leq C\left[\tau_1^\alpha + \left|\int_0^{\tau_1} b(t, X_t) \, dt\right| + \left|\int_0^{\tau_1} \gamma_t \, dW_t\right|\right]$$

$$\leq C[\tau_1^\alpha + \tau_1 + \tau_1^{1/3} \xi],$$

where $\xi \triangleq \sup_{0 \leq t \leq T} \frac{|\int_0^t \gamma_s \, dW_s|}{t^{1/3}}$. Then at least one of the following inequalities holds true:

$$\tau_1^\alpha \geq \frac{1}{6Cn}; \qquad \tau_1 \geq \frac{1}{6Cn}; \qquad \tau_1^{1/3} \xi \geq \frac{1}{6Cn}.$$

In any case, we have
$$\tau_1^{-1} \leq C\left[\frac{1}{T} + n^{1/\alpha} + n + n^3 \xi^3\right].$$

Since $|\gamma_t| \leq C$, applying Lemma 5.5, we get $E\{\xi^3\} \leq C < \infty$. Thus, $E\{\tau_1^{-1}\} \leq C_n < \infty$. Plugging this into (5.35) we get (5.34) and complete the proof for (i). □

DEPARTMENT OF MATHEMATICS
UNIVERSITY OF SOUTHERN CALIFORNIA
LOS ANGELES, CALIFORNIA 90089
USA
E-MAIL: jianfenz@usc.edu
URL: http://almaak.usc.edu/~jianfenz